*В данной работе рассматриваются матричные структуры арифметических процессоров, основанные на распределенной арифметике и многорядных кодах. Область применения - разработка суперкомпьютеров.*
***Ключевые слова и фразы****: распределенная арифметика; многорядные коды; сумматоры одноразрядных чисел; групповое суммирование; умножители; делители; накопительные сумматоры; нониусные преобразования; матричный арифметический процессор.*

*In this paper we consider the matrix structure of arithmetic processors based on distributed arithmetic in multi-row codes. Scope - development of supercomputers.*

**Щербаков В.И**., к.т.н. , x162@rambler.ru
Ph.D. **Shcherbakov V.I.,**
Институт проблем моделирования в энергетике НАН Украины.

## Сверхбыстрая арифметика в многорядных кодах
## Ultrafast a Distributed Arithmetic in multi-row codes

Непрерывно возрастающие требования к скорости обработки информации приводит к необходимости иметь огромные вычислительные ресурсы. Современные суперкомпьютеры, представляющие собой массив элементарных компьютеров, несмотря на распараллеливание алгоритмов, имеют большие временные потери на пересылки информации между компьютерами. Потребовались новые идеи в организации вычислительных процессов, в создании новых архитектур суперкомпьютеров, опирающиеся на достижения в технологии изготовления интегральных схем. В данной работе предлагается метод уменьшения временных затрат, связанных непосредственно с арифметическими процессами. Метод основывается на организации вычислений в многорядных кодах и использовании распределенной арифметики, позволяющей распараллеливать вычислительный процесс на уровне операций с отдельными разрядами.. Предлагается так же новая архитектура процессорных элементов, работающих в многорядных кодах, ориентированная на интегральную технологию изготовления микросхем

# 1. Форма представления чисел в многорядных кодах

В современной вычислительной технике используются в основном позиционные коды, в котором каждый разряд занимает определенную позицию в строке и имеет определенный вес. При сложении чисел, представленных в формате с плавающей запятой, необходимо выравнивать порядки чисел, в связи с чем теряются разряды мантисс, выходящие за пределы разрядной сетки, а, следовательно, и теряется информация. Возникают ошибки связанные с потерей точности результата при неполном пересечении мантисс чисел на числовой оси. Формат с фиксированной запятой позволяет складывать большие числа с малыми без потери информации, но имеет диапазон представления чисел, ограниченный разрядной сеткой. При увеличении разрядности кода увеличивается и время суммирования чисел, так как переносы необходимо распространять вплоть до самого старшего разряда. Увеличивать информационность кода можно путем представления операндов в формате матриц размером m x n, где m – число n- разрядных слов. Например, матрица из двух n- разрядных слов может быть выражена как :

| $a_{1,n}$ | … | $a_{1,i}$ | … | $a_{1,1}$ |
|---|---|---|---|---|
| $a_{2,n}$ | … | $a_{2,i}$ | … | $a_{2,1}$ |

Все коды, содержащиеся в определенном столбце такой матрицы, имеют определенный ( одинаковый ) вес. В дальнейшем, матрицы, составленные только из кодов, будем называть **кодовыми матрицами или многорядными кодами.** Так, выше приведенная матрица, в которой m=2, называться двухрядным кодом (2-код). Замечательным свойством такого формата является то, что при сложении матриц переносы не распространяются вдоль всей строки, а ограничиваются несколькими разрядами. Допустим числа А и В заданы в двухрядных кодах:

|  |  | | | | | | |
|---|---|---|---|---|---|---|---|
| A = | $a_{1,n}$ | $a_{1,n-1}$ | … | $a_{1,i}$ | … | $a_{2,1}$ | $a_{1,1}$ |
|  | $a_{2,n}$ | $a_{2,n-1}$ | … | $a_{2,i}$ | … | $a_{2,2}$ | $a_{2,1}$ |
|  | $b_{1,n}$ | $b_{1,n-1}$ | … | $b_{1,i}$ | … | $b_{2,1}$ | $b_{1,1}$ |
| B = | $b_{2,n}$ | $b_{2,n-1}$ | … | $b_{2,i}$ | … | $b_{2,2}$ | $b_{2,1}$ |

Тогда суммирование этих чисел будет представлять собой свертку (4→2) четырехрядного кода в двухрядный код [3]. Заметим, что сложение выполняется за фиксированное время, т.е. время выполнения сложения не зависит от числа разрядов. Обозначим сумму 4-х одноразрядных чисел $i$ –ого столбца в виде $C_i$. Если кодовые матрицы A и B представлены в двоичном формате, то $C_i$ будет представлена тремя разрядами соответствующего веса:

$C_i \rightarrow c_{i,3}\ c_{i,2}\ c_{i,1}$

Учитывая, что в пределах конкретного столбца положение разряда не имеет значения, то после первого этапа суммирования получим трехрядную кодовую матрицу частичных сумм **C** (трехрядный код):

| **C** = | - | - | $c_{n,1}$ | $c_{n-1,1}$ | ... | $c_{i,1}$ | ... | ... | $c_{2,1}$ | $c_{1,1}$ |
|---|---|---|---|---|---|---|---|---|---|---|
| | - | $c_{n,2}$ | $c_{n-1,2}$ | ... | $c_{i,2}$ | ... | ... | $c_{2,2}$ | $c_{1,2}$ | - |
| | $c_{n,3}$ | $c_{n-1,3}$ | ... | $c_{i,3}$ | ... | ... | $c_{2,3}$ | $c_{1,3}$ | - | - |

В результате свертки «3→2» получим двухрядную кодовую матрицу **S**, представляющую искомую сумму чисел A и B, в том же формате, что и входные операнды:

| **S** = | $s_{1,n+2}$ | $s_{1,n+1}$ | $s_{1,n}$ | ... | $s_{1,i}$ | ... | $s_{1,2} = c_{2,1}$ | $s_{1,1} = c_{1,1}$ |
|---|---|---|---|---|---|---|---|---|
| | $s_{2,n+2} = c_{n,3}$ | $s_{2,n+1}$ | $s_{2,n}$ | ... | $s_{2,i}$ | ... | $s_{2,2} = c_{1,2}$ | $s_{2,1} = 0$ |

Таким образом, суммирование двух операндов – матриц A и B размером 2 x n дает матрицу **S** размером 2 x (n+2). При этом независимо от разрядности n переносы распространяются только на три соседних разряда. Время суммирования определяется временем суммирования одного сумматора четырех одноразрядных чисел (в компрессоре 4→3) и одного сумматора трех одноразрядных чисел (в компрессоре 3→2). Число A, заданное в однорядном коде (например с фиксированной запятой), в двухрядном коде будет иметь следующий вид:

| A = | $a_{1,n}$ | ... | $a_{1,i}$ | ... | $a_{1,1}$ |
|---|---|---|---|---|---|
| | 0 | ... | 0 | ... | 0 |

Размер матриц-операндов может быть различным и ограничивается только технологическими возможностями изготовления интегральных схем [1]. Если старшие разряды строк в матрице m x n (старший столбец) отвести под знаковые разряды, то такая кодовая матрица может представлять числа в дополнительном коде.  Следовательно, операции вычитания в многорядных кодах (так же как и операции суммирования) сводятся к свертке матрицы, выраженной в дополнительных кодах. При обработке больших массивов данных формирование дополнительных кодов требует больших временных затрат.  Кроме того, дополнительный код априори  требует задавать определенный размер разрядной сетки, что при обработки больших массивов чисел возникают проблемы с переполнением разрядной сетки. Для таких задач мы предлагаем четырехрядные операнды-матрицы, в которых положительные двухрядные числа позиционируются в двух верхних (или нижних) рядах, а отрицательные в двух нижних (или верхних) рядах.

$$A = \begin{array}{c} A_+ = \\ \\ A_{(-)} = \end{array} \begin{array}{|c|c|c|c|c|c|c|} \hline a_{1,n} & a_{1,n-1} & \ldots & a_{1,i} & \ldots & a_{2,1} & a_{1,1} \\ \hline a_{2,n} & a_{2,n-1} & \ldots & a_{2,i} & \ldots & a_{2,2} & a_{2,1} \\ \hline a_{1,n} & a_{1,n-1} & \ldots & a_{1,i} & \ldots & a_{2,1} & a_{1,1} \\ \hline a_{2,n} & a_{2,n-1} & \ldots & a_{2,i} & \ldots & a_{2,2} & a_{2,1} \\ \hline \end{array}$$

В таком формате при суммировании положительные числа складываются с положительными числами, а отрицательные – с отрицательными.  Результат сложения так же представляется в таком же формате.

$$S = A + B = \begin{array}{|c|} \hline A_{(+)} \\ \hline A_{(-)} \\ \hline \end{array} + \begin{array}{|c|} \hline B_{(+)} \\ \hline B_{(-)} \\ \hline \end{array} = \begin{array}{|c|} \hline S_{(+)} = A_{(+)} + B_{(+)} \\ \hline S_{(-)} = A_{(-)} + B_{(-)} \\ \hline \end{array}$$

При вычитании  необходимо в вычитаемом поменять местами матрицы (верхний код с нижним кодом), после чего произвести сложение:

$$S = A - B = \begin{array}{|c|} \hline A_{(+)} \\ \hline A_{(-)} \\ \hline \end{array} + \begin{array}{|c|} \hline B_{(-)} \\ \hline B_{(+)} \\ \hline \end{array} = \begin{array}{|c|} \hline S_{(+)} = A_{(+)} + B_{(-)} \\ \hline S_{(-)} = A_{(-)} + B_{(+)} \\ \hline \end{array}$$

Очевидно, что при вычитании в данном случае не требуются дополнительные коды. Все коды A(+), A(-), B(+), B(-), S(+), S(-) в данном формате не имеют знаковых разрядов, а их форматы есть прямые коды. При этом вес разрядов в  младшем столбце матрицы назначается специальным кодом (сопутствующим матрице), значение которого кратно байту. Это позволяет позицию запятой, которая  отделяет целую часть числа от дробной части, определять виртуально по значению веса разрядов младшего столбца.  Этот формат  можно назвать как  «формат с виртуальной позицией точки».   Например, если вес разрядов в младшем столбце задан как $2^{-128}$, то запятая будет располагаться на 127 разрядов (столбцов) влево от младшего столбца матрицы.  Очевидно, что определив вес разрядов в младшем столбце, все операции должны производится с матрицами-операндами имеющими те же веса разрядов в младшем столбце. При таком представлении операндов-матриц операцию суммирования с накоплением можно производить без потери точности, задав вес разрядов в младшем столбце равным минимально возможному значению входного числа. Магистрали для входных и выходных данных при таком формате должны представлять собой 4-х шинные магистрали. Все арифметические операции (суммирование, вычитание, умножение, деление), в конечном итоге, сводятся к одному типу операции - операции сложения, а в предлагаемом формате - к операции свертки многорядного кода. Рассмотрим вычисление оператора группового суммирования $S = \Sigma A_i$ и оператора суммы парных произведений $F = \Sigma A_i B_i$ (i=1…m). Эти операторы составляют основу вычисления разностных уравнений, уравнений типа «бабочка» в спектральном анализе, матричных уравнениях, в уравнениях для цифровых фильтров и т.д. Скорость вычисления этих операторов во многом определяют скорость всего вычислительного процесса при обработке потока данных.

## 2. Вычисление оператора группового суммирования

Оператор группового суммирования

$$S = \sum_{i=0}^{m-1} A_i \quad , \quad A_i = \sum_{j=0}^{n-1} a_{i,j} g^j \qquad (2.1)$$

может быть выражен в матичном виде как:

$$S = A\,Q$$

или

$$S = \begin{array}{|c|c|c|c|c|} \hline a_{0,n-1} & \ldots & a_{0,j} & \ldots & a_{0,0} \\ \hline \ldots & \ldots & \ldots & \ldots & \ldots \\ \hline a_{i,n-1} & \ldots & a_{i,j} & \ldots & a_{i,0} \\ \hline \ldots & \ldots & \ldots & \ldots & \ldots \\ \hline a_{m-1,n-1} & \ldots & a_{m-1,j} & \ldots & a_{m-1,0} \\ \hline \end{array} \begin{array}{|c|} \hline g^{n-1} \\ \hline \ldots \\ \hline g^j \\ \hline \ldots \\ \hline g^0 \\ \hline \end{array}$$

or

$$S = \begin{array}{|c|c|c|c|c|} \hline a_{0,n-1} & \ldots & a_{0,j} & \ldots & a_{0,0} \\ \hline \ldots & \ldots & \ldots & \ldots & \ldots \\ \hline a_{i,n-1} & \ldots & a_{i,j} & \ldots & a_{i,0} \\ \hline \ldots & \ldots & \ldots & \ldots & \ldots \\ \hline a_{m-1,n-1} & \ldots & a_{m-1,j} & \ldots & a_{m-1,0} \\ \hline \end{array} \quad \begin{array}{|c|} \hline g^{n-1} \\ \hline \ldots \\ \hline g^{j} \\ \hline \ldots \\ \hline g^{0} \\ \hline \end{array}$$

Матрица **A**, сформированная из m однорядных кодов, называется кодовой матрицей или m-рядным кодом.

Выразим оператор (2.1) в избыточной системе счисления как:

$$S = \sum_{i=0}^{m-1}\sum_{j=0}^{n-1} a_{i,j} g^{j} = \sum_{j=0}^{n-1} S_j g^{j} \tag{2.2}$$

где $S_j = \sum_{j=0}^{m-1} a_{i,j}$ - сумма одноразрядных чисел j-ого столбца матрицы **A**

Поскольку m-рядный код S необходимо свернуть в двухрядный код, то процесс суммирования будет происходить в несколько этапов (слоев). Сумму S, так же как и кодовые матрицы, соответствующие этапу свертки, будем обозначать индексом в круглых скобках. Так, кодовая исходная матрица **A** в выражении (2.1) будет иметь обозначение $\mathbf{A}_{(1)}$, а сумма одноразрядных чисел в j-ом столбце этой матицы будет иметь обозначение $S_{(1),j}$. Аналогично для размерностей матрицы $\mathbf{A}_{(1)}$: $m \to m_{(1)}$ ; $n \to n_{(1)}$ ….

Выразим $S_{(1),j}$ в канонической позиционной системе счисления с основанием q:

$$S_{(1),j} = \sum_{h=0}^{m_{(2)}-1} a_{i,j} q^{h} \tag{2.3}$$

где $\quad m_{(2)} = ]\log_q (1+S_{(1),j}^{max})[ \; = \; ]\log_q (1+m_{(1)}q - m_{(1)})[ \quad ,$ (2.4)

$S_{(1),j}^{max} = m_{(1)}(q-1)$ - максимальное значение суммы одноразрядных чисел в j-ом столбце; ]...[ - знак округления до ближайшего большего целого

Запишем выражение (2.2) с учетом (2.3) в матричном виде

$$S = A_{(2)} Q,$$

где кодовая матрица $A_{(2)}$ имеет следующий вид:

$A_{(2)} =$

| | | | | | $a_{0,n-1}$ | … | $a_{0,j}$ | … | $a_{0,1}$ | $a_{0,0}$ |
|---|---|---|---|---|---|---|---|---|---|---|
| | | | | $a_{1,n-1}$ | … | $a_{1,j}$ | … | $a_{1,1}$ | $a_{1,0}$ | |
| | | | … | … | … | … | … | … | | |
| | | $a_{i,n-1}$ | … | $a_{i,j}$ | … | $a_{i,1}$ | $a_{i,0}$ | | | |
| | … | … | … | … | … | … | | | | |
| $a_{m-1,n-1}$ | … | $a_{m-1,j}$ | … | $a_{m-1,1}$ | $a_{m-1,0}$ | | | | | |

В матрице $A_{(2)}$: $n = m_{(2)}$; $m = n_{(1)}$

Поскольку матрица $A_{(2)}$, также как и матрица $A_{(1)}$, отражают сумму S, то расположение разрядов $a_{i,j}$ в пределах j-ого столбца не имеет значения. Поэтому, сместив соответствующие разряды в пределах собственных столбцов в матрице, сформируем матрицу $A_{(2)}$ в форме трапеции:

$A_{(2)} =$

| $a_{m-1,n-1}$ | … | $a_{i,n-1}$ | … | $a_{1,n-1}$ | $a_{0,n-1}$ | … | $a_{0,j}$ | … | $a_{0,1}$ | $a_{0,0}$ |
|---|---|---|---|---|---|---|---|---|---|---|
| | … | … | … | … | … | … | … | … | … | |
| | | … | … | … | … | … | … | … | | |
| | | | … | … | … | … | … | | | |
| | | | | $a_{m-1,j}$ | … | $a_{m-1,0}$ | | | | |

В этой матрице: $m = m_{(2)}$; $n = n_{(1)} + m_{(2)} - 1$.

Разрядность числа, находящегося в (m-1)–ой строке трапецевидной матрицы $A_{(2)}$ минимальна и равна $n_{min} = n_{(1)} - m_{(1)} + 1$.

Матрица $A_{(2)}$ представляет $m_{(2)}$–рядный код, который соотносится с $m_{(1)}$–рядным кодом логарифмической зависимостью (2.4). Особенностью кодовой матрицы $A_{(2)}$ является то, что она также отражает матрицу частичных произведений двух однорядных кодов.

Например, произведение $F = BC$, где $B = b_{n-1}q^{n-1} + … b_i q^i + … b_0 q^0$, $C = c_{n-1}q^{n-1} + … c_i q^i + … c_0 q^0$, (i-целые положительные числа) можно представить в форме кодовой матрицы F следующего вида:

$\mathbf{F}_{(2)} =$

|  |  |  |  | $b_0 \& c_{n-1}$ | ... | $b_0 \& c_i$ | ... | $b_0 \& c_0$ |
|---|---|---|---|---|---|---|---|---|
|  |  |  | ... | ... | ... | ... | ... |  |
|  |  | ... | ... | ... | ... | ... |  |  |
|  | ... | ... | ... | ... | ... |  |  |  |
| $b_{n-1} \& c_{n-1}$ | ... | $b_{n-1} \& c_i$ | ... | $b_{n-1} \& c_0$ |  |  |  |  |

В связи с этим все последующие этапы свертки многорядного кода будут отражать также процесс умножения двух однорядных кодов. Каждая последующая кодовая матрица $\mathbf{A}_{i+1}$ в процессе свертки образуется из матрицы $A_i$ при операции суммирования одноразрядных чисел (2.3), при этом рядность (число строк) каждой последующей матрицы связана с рядностью предыдущей логарифмической зависимостью (2.4).

В таблице 2.1 приведены значения рядности $m_i$ кодовых матриц, образующихся на каждом этапе, при свертке исходного m-рядного кода (для двоичной системе счисления).

Таблица 2.1

| $m_{(1)}$ – число строк | $m_{(2)}$ | $m_{(3)}$ | $m_{(4)}$ | L - число этапов |
|---|---|---|---|---|
| $2^2$ -1 | 2 | - | - | 1 |
| $2^2$... $(2^3$ -1) | 3 | 2 | - | 2 |
| $2^3$... $(2^4$ -1) | 4 | 3 | 2 | 3 |
| $2^4$... $(2^5$ -1) | 5 | 3 | 2 | 3 |
| $2^5$... $(2^6$ -1) | 6 | 3 | 2 | 3 |
| $2^6$... $(2^7$ -1) | 7 | 3 | 2 | 3 |

Анализ данных, приведенных в таблице 2.1, показывает, что функция L=f($m_i$) ступенчата, где L- число этапов свертки исходной кодовой матрицы в двухрядный код. Значение L в определенном диапазоне не зависит от числа слагаемых. Например, 64 или 127 слагаемых можно свернуть в двухрядный код за 3 этапа. Этот факт (при конструировании аппаратуры) указывает на возможность ускорить вычислительный процесс путем выбора оптимального размера исходной кодовой матрицы. .

Аппаратурные затраты на свертку многорядных кодов в основном определяются количеством и числом входов сумматоров одноразрядных чисел (дальше по тексту - ОСЧ). Рассчитать эти затраты можно на основе следующих соотношений:

для матрицы $\mathbf{A}_{(2)}$ :

$n_{(2),i} = n_{(2)}^{max} - \Delta n_{(2),i}$ - разрядность i-ой строки

$\Delta n_{(2),i} = 0$ , если  i=0;

$\Delta n_{(2),i} = 2^i$ , если $i < m_{(2)}$,  i= 1, …$(m_{(2)} – 1)$ – номер строки

$N_{(2),v} = 2$,  если $v < m_{(2)}$,   - число столбцов, имеющих длину v

$N_{(2),v} = n_{(2)}^{min}$ , если $v = m_{(2)}$,  $v = (1, … m_{(2)} )$ – длина столбца (число разрядов)

$n_{(2)}^{max} = n_{(1)} + n_{(2)} -1$

$n_{(2)}^{min} = n_{(1)} - m_{(2)} +1$

для матрицы $\mathbf{A}_{(3)}$ :

$n_{(3),i} = n_{(3)}^{max} - \Delta n_{(3),i}$   - разрядность i-ой строки

$\Delta n_{(3),i} = 0$ ,  если  i=0;
$\Delta n_{(3),i} = 2^i$ , если $i < m_{(3)}$,  i= 1, …$(m_{(3)} – 1)$ – номер строки
$N_{(3),v} = 2$,  если $v < m_{(3)}$,    - число столбцов, имеющих длину v

$N_{(3),v} = n_{(3)}^{min}$ , если $v = m_{(3)}$,   $v = (1, … m_{(3)} )$ – длина столбца (число разрядов)
$n_{(3)}^{max} = n_{(2)}^{max}$
$n_{(3)}^{min} = n_{(3)}^{max} - q^{m(3)} +1$

Временные затраты на вычисление оператора группового суммирования в основном определяются  быстродействием сумматора одноразрядных чисел (ОСЧ) и числом этапов свертки. В свою очередь, быстродействие ОСЧ определяется числом слагаемых.  Очевидно, что наибольшее число слагаемых ($m_{(1)}$) в одном столбце содержит матрица $\mathbf{A}_{(1)}$.
Время свертки исходной матрицы размером  n x $m_1$ в двухрядный код определяется числом этапов свертки L, и временем $t_m$ суммирования в ОСЧ в каждой конкретной матрице, которое, зависит от числа слагаемых m.     Если, например, число слагаемых (при любом n) в исходной матрице не превышает 63 ($m_1$=63), то время свертки такой матрицы в двухрядный код с помощью табличных ОСЧ  будет  равно:

$$t_{mxn} = t_{m1} + t_{m2} + t_{m3} + t_{cc} = 13t_\& \qquad (2.5)$$
$$t_{m1} = 6t_\& , \; t_{m2} = 3t_\&, \; t_{m3} = t_\& , \; t_{cc} = 3t_\&$$

где $t_\&$ - задержка в одном элементе конъюнкции,

$t_{cc}$ - задержка в шифраторе сумматора

$m_1=63$, $m_2=6$, $m_3=3$, $m_4=2$

## 3. Табличные сумматоры одноразрядных чисел

Из анализа кодовых матриц $A_{(1)}$, $A_{(2)}$, $A_{(3)}$ видно, что минимальное число входов в одноразрядном сумматоре равно 3, а максимальное равно числу слагаемый в исходной матрице ($m_{(1)}$). Если $m_{(1)}$ велико, то реализовать одноразрядный сумматор с табличной структурой (в пределах допустимых аппаратурных затрат) можно в форме древовидной структуры. Ниже будут рассматриваться только табличные сумматоры, поскольку их структура наиболее полно отвечает требованиям интегральной технологии. Рассмотрим элементарный табличный сумматор, имеющий три входа и два выхода (компрессор 3→2)

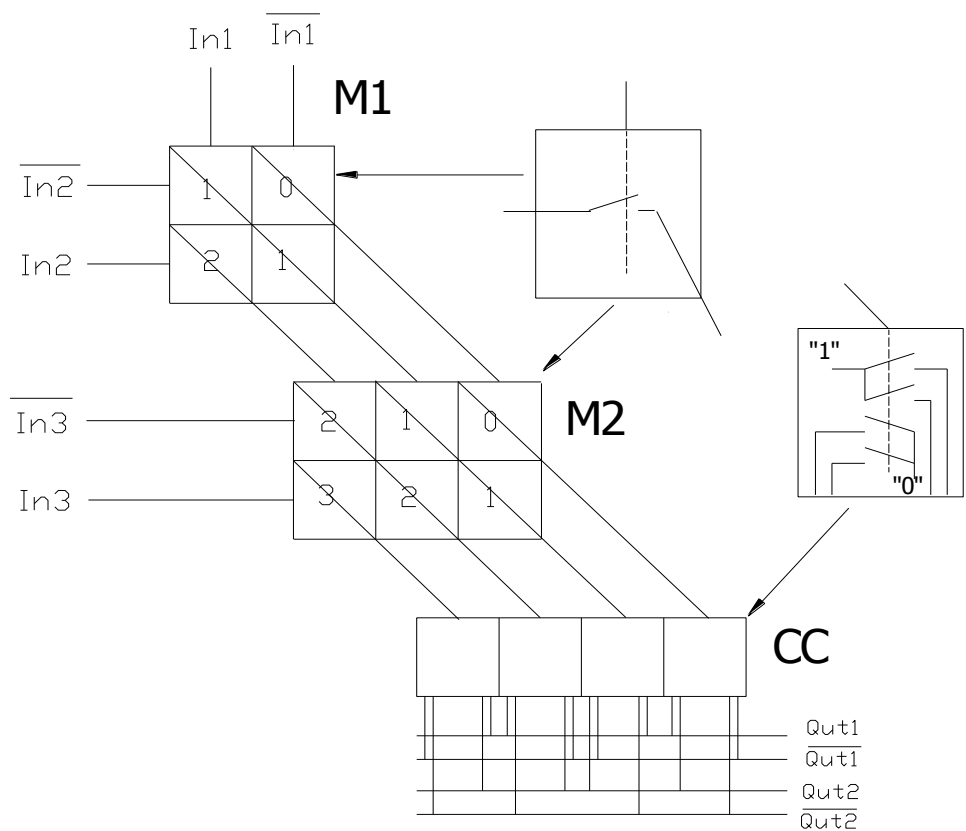

Fig.1

На рис. 1 показана структура компрессора 3→2 . Элементы таблиц M1,M2 представляют собой двухвходовые элементы конъюнкции (AND). Элементы шифратора СС можно выполнить на транзисторных элементах типа «ключ». Заметим, что элементы таблиц М2, М3 также могут быть выполнены в форме ключа. На первый вход ключа (например, затвор) подается разрешающий сигнал, а на второй вход (например, сток или исток).подается информационный сигнал. Время сложения в данном компрессоре равно $t_{3-2} = 2 t_\& + t_k$ , где $t_\&$, $t_k$ время включения двухвходового элемента AND и одного ключа. Если весь сумматор построен на ключах то: $t_{3-2} = 3t_k$ . При этом аппаратурные затраты составляют 26 ключей.

Рассмотрим сумматор одноразрядных чисел, имеющий m – входов (компрессор m-2) . Такой компрессор будет иметь древовидную структуру.

Для упрощения пояснений будем рассматривать только квадратные таблицы, составляющие древовидную структуру сумматора, показанную на рис.2

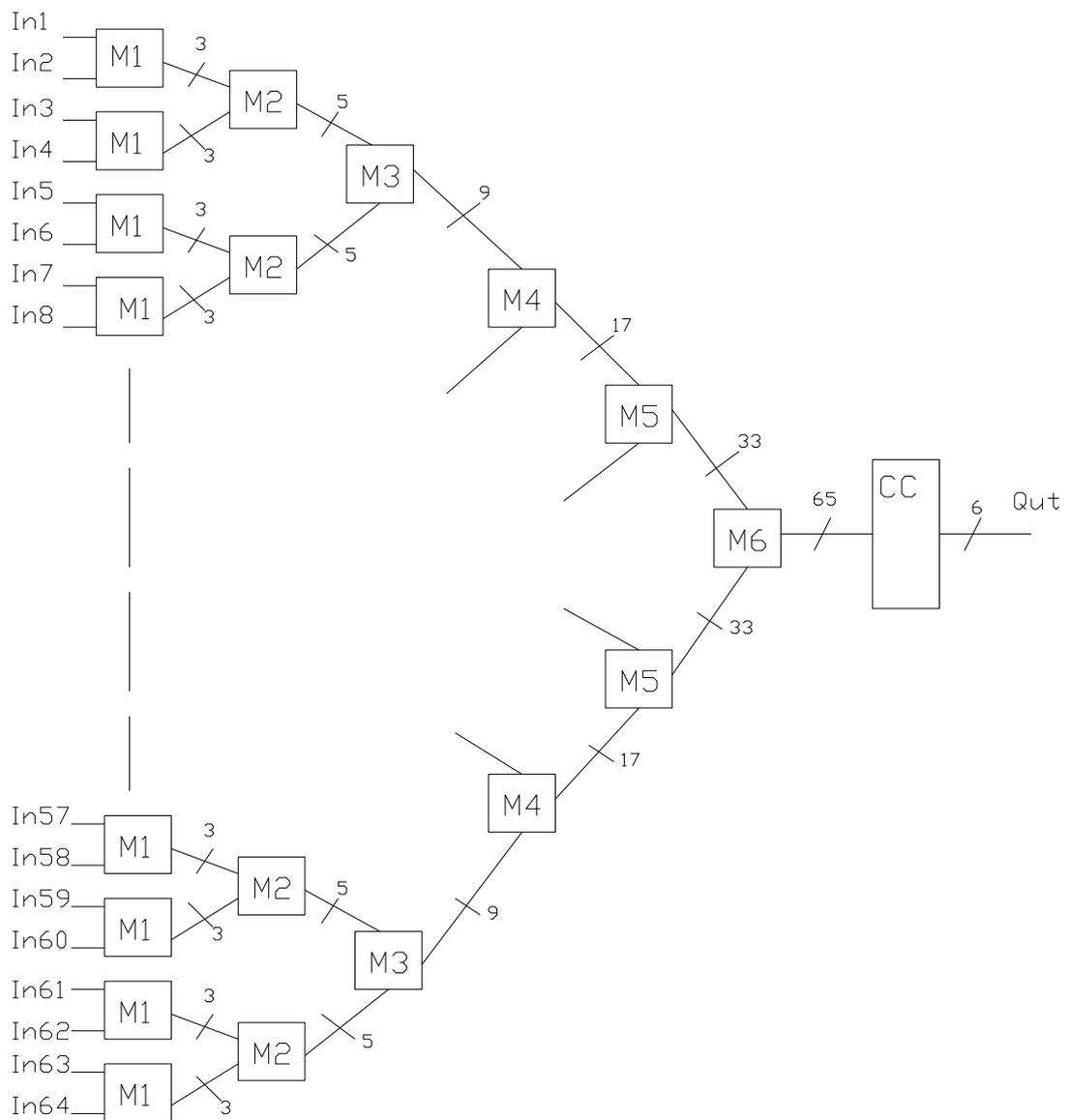

Fig.2

Древовидная структура сумматора для 64 одноразрядных чисел содержит 6 типов таблиц M1…M6.   В общем случае (для квадратных таблиц) количество типов таблиц N связано с числом слагаемых m соотношением: $m = 2^N$. При этом таблицы имеют следующие размеры : M1 = 2 x 2, M2 = 3 x 3 , M3 = 5 x 5 , M4 = 9x9. M5 = 17 x 17, M6 = 33 x 33. Эти размеры (для квадратных таблиц) не зависят от числа слагаемых .   Количество таблиц $V_i$, где i номер типа таблицы, связано с числом (m) слагаемых (кратным $2^i$) следующим соотношением: $V_i = m/2^i$.  Так, для m=64  имеем: M1= 32; M2 = 16; M3 = 8; M4 = 4; M5 = 2; M6 = 1. Поскольку таблица, для любой размерности,  имеет временную задержку $t_\&$, равную  задержке в одном двухвходовом  элементе конъюнкции, то общее время сложения 64 одноразрядных чисел будет равно   $t_s = 7 t_\&$ , с учетом времени задержки в шифраторе СС равной $t_{cc}= t_\&$.   В диапазоне входных сигналов (3-128) время сложения m – одноразрядных чисел можно определить из таблицы 3.1

Таблица 3.1

| Число входов | 3-4 | 5-7 | 8-16 | 17-32 | 33-64 | 65-128 |
|---|---|---|---|---|---|---|
| Время суммирования | $2t_\& + t_{cc}$ | $3t_\& + t_{cc}$ | $4t_\& + t_{cc}$ | $5t_\& + t_{cc}$ | $6t_\& + t_{cc}$ | $7t_\& + t_{cc}$ |

Из данной таблицы следует, что время сложения выражается ступенчатой функцией, что следует учитывать при конструировании устройств свертки многорядных кодов.   Аппаратурные затраты табличных сумматоров одноразрядных чисел для диапазона входов (3-32)  сведены в таблицу 3.2.

Таблица 3.2

| M, $\sum_\&$ | Число входов сумматора (m). | | | | | | | | | | | | | | |
|---|---|---|---|---|---|---|---|---|---|---|---|---|---|---|---|
|  | 3 | 4 | 6 | 8 | 10 | 12 | 14 | 16 | 18 | 20 | 22 | 24 | 26 | 28 | 30 | 32 |
| M1 | 1 | 2 | 3 | 4 | 5 | 6 | 7 | 8 | 9 | 10 | 11 | 12 | 13 | 14 | 15 | 16 |
| M2 | 1 | 1 | 1 | 2 | 2 | 3 | 3 | 4 | 4 | 5 | 5 | 6 | 6 | 7 | 7 | 8 |
| M3 | - | - | - | - | 1 | 1 | 1 | 2 | 2 | 2 | 3 | 3 | 3 | 3 | 4 | 4 |
| M4 | - | - | - | - | 1 | 1 | 1 | 1 | 1 | 1 | 1 | 1 | 2 | 2 | 2 | 2 |
| M5 | - | - | - | - | - | - | - | 1 | 1 | 1 | 1 | 1 | 1 | 1 | 1 | 1 |
| $\sum_\&$ | 10 | 17 | 36 | 59 | 90 | 121 | 158 | 199 | 254 | 301 | 354 | 411 | 476 | 541 | 582 | 687 |

В данной таблице указано число матриц типа M1-M5 и суммарное число двухвходовых элементов AND ($\sum_\&$), входящих в шифраторы СС, для конкретного  числа входов (m) сумматора одноразрядных чисел. Для 64–х одноразрядных чисел число двухвходовых элементов AND,  входящих в шифратор С С, будет равно $\sum_\& = 2463$ .

Свертка многрядных кодов  на основе  компрессоров 6-2 и 4-2  описана в работах  [2, 5].

## 4. Операция группового суммирования с накоплением

Запишем операцию группового суммирования с накоплением в следующем виде:

$$S = \sum_{i=1}^{w} S_i, \qquad S_i = \mathbf{A}_i + S_{i-1} \qquad (4.1)$$

где $S_{i-1}$ – текущее значение суммы; $\mathbf{A}_i$ – входная матрица-операнд.

В данном случае на каждой итерации на вход поступают матрица-операнд $\mathbf{A}^i$, размером (n x m), и результат предыдущей свертки многорядного кода $S_{i-1}$. Учитывая, что (при свертке многорядного кода) двухрядному коду всегда предшествует трехрядный код, то в выражении (4.1) $S_{i-1}$ можно представлять трехрядным кодом. Это позволяет сократить число этапов свертки на один этап и ускорить процесс свертки. На (w-1) – итерации необходимо свернуть трехрядный код $S_{w-1}$ в двухрядный код. Рассмотрим наиболее простую архитектуру накапливающего сумматора (рис.3), в котором операнд $A_i$ представлен однорядным кодом, а текущая сумма $S_{i-1}$ представлена двухрядным кодом.

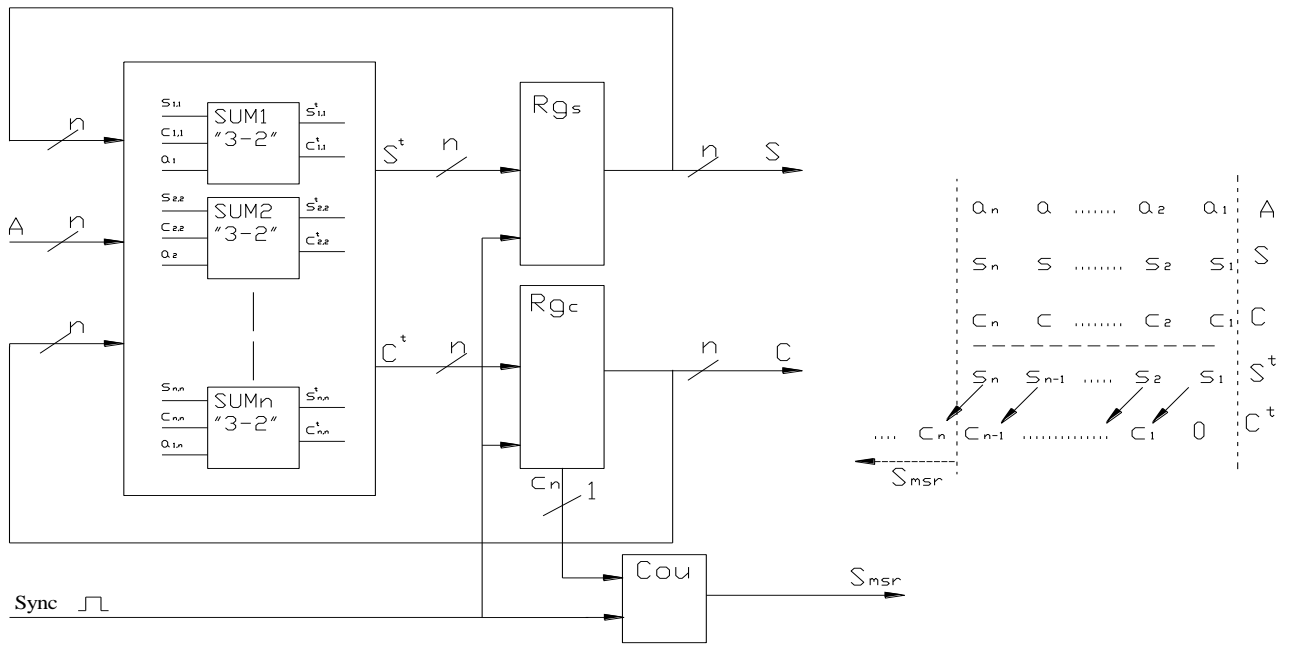

Fig. 3

Предполагается отдельное суммирование положительных и отрицательных чисел. Поэтому все числа представлены в прямых кодах без знаковых разрядов. Очевидно, что в данном случае основу такого сумматора будет составлять линейка из n – сумматоров одноразрядных чисел типа $S_{3 \to 2}$, роль которого может выполнять полный сумматор.

В структуру такого сумматора входят так же два регистра $Rg_s$, $Rg_c$ (построенных триггерах типа "мастер-помощник"). Регистр $Rg_s$ принимает от сумматоров разряды суммы, а регистр $Rg_c$ - разряды переносов. Расположение разрядов при суммировании показано на этом же рисунке в форме отдельной выноски. Время суммирования одного n-разрядного числа $A_i$ определяется временем суммирования в одном сумматоре $S_{3 \to 2}$ и временем записи информации в регистр.

Поскольку переносы при суммировании в данной структуре распространяются только на один разряд, то все входные и выходные коды можно разметить в n- разрядной сетке. Переполнение при сложении определяется переносом $c_n$ n-го (старшего) сумматора $S_{3 \to 2}$. Этот перенос поступает на вход двоичного счетчика импульсов Cou с параллельным двоичным выходом. В такой структуре накапливающий сумматор может иметь сколь угодно большой диапазон выходного кода, и суммирование производится без потери точности вычислений. При этом если, например, требуется накапливать малоразрядные числа $A_i$ с весьма большим числом итераций (циклов накопления w), то такая структура накопителя будет иметь минимальные аппаратурные затраты, поскольку линейка одноразрядных сумматоров и регистров так же будет малоразрядной. Если операнд A представлен двухрядным кодом, то потребуется свертка четырехрядного кода в двухрядный код в два этапа, при этом переносы будут распространяться на два разряда. В этом случае разрядность регистра $R_{gs}$ необходимо увеличить на один разряд, выходой разряд которого $s_{n+1}$ должен суммироваться по модулю два с выходным разрядом $c_{n+1}$ регистра $R_{gc}$, а результат суммирования поступает на вход счетчика Cou.

## 5. Оператор суммы парных произведений

Прежде чем рассматривать процесс вычисления оператора суммы парных произведений

$$F = \sum_{i=1}^{m} A_i B_i \qquad (5.1)$$

рассмотрим произведение P пары чисел A, B, заданных в дополнительном коде :

$$P = A B ,$$

где $A = -a_{(0)}q^0 + \sum_{i=1}^{n} a_i q^{-i}$, $B = -b_{(0)}q^0 + \sum_{j=1}^{n} b_j q^{-j}$,

$a_{(0)}$, $b_{(0)}$ – знаковые разряды, а числа A, B могут быть как нормализованные, так и ненормализованные.

Запишем произведение этих чисел в следующем виде:

$$P = AB = D + E - F - C, \qquad (5.2)$$

где $D = a_{(0)}b_{(0)}q^{(0)}$; $E = \sum_{i=1}^{n}\sum_{j=1}^{n} a_i b_j q^{-(i+j)}$; $F = \sum_{j=1}^{n} a_0 b_j q^{-j}$; $C = \sum_{i=1}^{n} b_0 a_i q^{-i}$

Слагаемое E представляет собой матрицу частичных произведений, не зависящую от знаковых разрядов, а слагаемые D,F,C – корректирующие числа для получения произведения в дополнительном коде. При $a_0 = 0$ имеем D=F=0, а при $b_0 = 0$ получим D=C=0. Заменим операцию вычитания в выражении (5.2) операцией сложения, для этого представим числа F и C в дополнительном коде:

$$P = AB = D + E + F_{ac} + C_{ac} \qquad (5.3)$$

Кодовая матрица P, сформированная в соответствии с выражением (5.3)., имеет следующий вид:

| $2^0$ | $2^{-1}$ | $2^{-2}$ | $2^{-3}$ | ... | ... | ... | ... | ... | ... | $2^{-2n}$ |
|---|---|---|---|---|---|---|---|---|---|---|
| $b_0a_0$ | **$b_1a_0$** | $b_1a_1$ | $b_2a_1$ | $b_3a_1$ | ... | ... | ... | ... | ... | $b_na_n$ |
| | **$b_0a_1$** | **$b_2a_0$** | $b_1a_2$ | $b_2a_2$ | $b_2a_3$ | ... | ... | ... | ... | |
| | | **$b_0a_2$** | **$b_3a_0$** | ... | ... | ... | ... | ... | | |
| | | | **$b_0a_3$** | ... | ... | ... | ... | | | |
| | | | | ... | **$b_n a_0$** | ... | | | | |
| | | | | 1 | **$b_0a_n$** | | | | | |

Жирным шрифтом в данной матрице показаны инвертированные конъюнкции. Данная матрица P, представляющая собой матрицу частичных произведений, соответствует матрице $A_{(2)}$, рассмотренной в разделе 2. Поэтому процесс свертки данной матрицы в двухрядный код и оценка временных и аппаратурных затрат соответствует свертке матрицы

А(2), рассмотренной выше. Число этапов L свертки матрицы Р можно определить из табл.2.1. Заметим, что если операнды А и В заданы в прямом коде без знаковых разрядов. При этом в приведенной выше матрице частичных произведений необходимо исключить единицу в последней строке. При этом инвертированные разряды заменить на не инвертированные, при этом знаковые разряды $a_0$, $b_0$ считать информационными с весом $2^0$. Время умножения определяется временем формирования матрицы Р, равной задержке $t_\&$ в одном элементе AND, и временем свертки m-рядного кода матрицы Р в двухрядный код. Например, время умножения двух 63-х разрядных кодов (см.2.5 ) будет равно 14 $t_\&$ .

Представим оператор суммы парных произведений (5.1) в рекуррентной форме:
$$F_i = F_{i-1} + A_i B_i = F_{i-1} + \mathbf{P}_i$$

В связи с тем, что функция $L=f(m_i)$ ступенчатая, то операнд $F_{i-1}$ можно включить в кодовую матрицу Р. Это позволяет исключить операцию сложения в данном выражении. Чтобы не увеличивать число этапов свертки (соответственно и время свертки) матрицы Р, необходимо определить этап, в который включается $F_{i-1}$. Заметим, что кодовую матрицу Р можно свернуть в двухрядный код с помощью табличных малоразрядных умножителей ( на основе малоразрядных ПЗУ) [1] или компрессоров 6→2; 4→2; 3→2. [2,3,5].

## 6. Операция деления

Операцию деления запишем в следующем виде:
$$y = x/z ,$$
или $\quad y^t = x$ , где $y^t = y\, z$ \hfill (6.1)

Из выражения (6.1) следует, что в операции деления участвуют операции умножения и сравнения. В современных алгоритмах, реализующих операцию деления, используется способ вычисления «цифра за цифрой», в котором на каждой итерации определяется один разряд частного. Рассмотрим ускоренный способ деления, в котором на каждой итерации вычисляются k- разрядов частного. При этом величина k ограничивается только технологическими возможностями интегральной технологии. Суть способа состоит в векторной операции сравнения x с цифровой шкалой $\mathbf{Y}_1^{\,t}$ на первой итерации и вычисления «невязки», после чего в операции сравнения на второй итерации участвует вычисленная «невязка» и т.д. Рассмотрим этот процесс более подробно. Введем понятие цифровой линейной шкалы. Под цифровой шкалой **Y** будем понимать линейную шкалу, имеющую $q^n$ – цифровых делений (включая нулевое деление), где q-основание системы счисления. Все весовые значения делений на такой шкале представлены n- разрядными кодами. Значения каждых двух соседних делений на такой шкале отличаются друг от друга на один «квант», т.е. на величину младшего разряда кода. Таким образом, на такой шкале будут присутствовать все возможные значения n-разрядного числа и , следовательно, определение неизвестного n-разрядного числа x путем сравнения его с такой шкалой может производится за одну итерацию.

Поскольку какая шкала должна храниться в памяти или формироваться в процессе вычислений, то при большой разрядности числа аппаратурные затраты на такую шкалу будут недопустимо большими. Вместо одной «длинной» цифровой шкалы **Y** , имеющей $q^n$ делений: **Y** → $y_1$ $y_2$ … $y_i$ … $y_n$ , в операциях сравнения используем такой набор цифровых шкал, в которых «длина» каждой последующей шкалы перекрывает интервал между двумя соседними делениями предыдущей шкалы.

**Y** → $(y_0\ y_1\ ....\ y_k)_1\ (y_{k+1}\ y_{k+2}\ ....\ y_{k+2k})_2\ ....\ (y_{m-k}\ ....\ y_n)_m$ → $Y_1\ Y_2…Y_j…Y_m$

Таким образом , если сформировать только первую шкалу $Y_1$, то все остальные нониусные шкалы можно получать путем сдвига кодов делений шкалы на k-разрядов вправо в начале каждой последующей итерации, т.е.

$$\mathbf{Y} = \sum q^{-k(j-1)} \mathbf{Y}_1, \text{ где } j = 1, 2….m \qquad (6.2)$$

Из выражения (6.1) следует, что в операции сравнения участвует не цифровая шкала **Y**, а ее модификация $\mathbf{Y}^t$ , представляющая собой произведение скалярной величины z на вектор **Y.** То есть все деления шкалы **Y** необходимо умножить на n-разрядный код числа z , при этом значения делений шкал **Y** и $\mathbf{Y}^t$ однозначно связаны (таблично) соответствующими номерами делений. В результате сравнения скалярной величины x с вектором $\mathbf{Y}_1^t$ (6.1) на первой итерации получим вектор разностей $\Delta$ . Анализируя знаковые разряды $_s\Delta$ , получим унитарный код, который однозначно определяет номер деления на шкале $\mathbf{Y}_1^t$ и , соответственно, на шкале **Y**. На выходе дешифратора единичного унитарного кода будут сформированы первые старшие k разрядов частного:

$x\ (><=)\ \mathbf{Y_1^t} \to \Delta_1 \to \Delta_{1,h} \to y^t_{1,h} \to y_{1,h} \to (y_0\ y_1\ ....\ y_k)_1$

где индекс h обозначает номер ближайшего деления (по недостатку), вес которого примерно равен x. Одновременно с дешифрацией унитарного кода определяется и величина невязки (по недостатку):

$$\Delta_1 = y^t_{1,h} - x$$

На второй итерации в операции сравнения со шкалой $\mathbf{Y}_2^t$ участвует невязка $\Delta_1$. Вместо формирования шкалы $\mathbf{Y}_2^t$ путем сдвига кодов ее делений вправо на k разрядов (6.2) можно сдвигать на k разрядов влево невязку $\Delta_1$. Математическую модель процесса деления по невязкам представим в следующем виде:

$$\Delta_{j-1,h} \; (><=) \; \mathbf{Y_1^t} \to \Delta_j \to \Delta_{j,h} \to y^t_{j,h} \to y_{j,h} \to (y_0 \; y_1 \; \ldots \; y_k)_j$$

$$\Delta_{j,h} = y^t_{j-1,h} - \Delta_{j-1} \quad \text{для} \quad j = 1\ldots m$$

$$\Delta_{j+1,h} = q^k \Delta_{j,h} \quad \text{для} \quad j > 1$$

$$\Delta_{j-1,h} = x \quad \text{для} \quad j = 1$$

Операции умножения и сравнения можно выполнять одновременно, используя свертку многорядного кода. При этом вместо операции вычитания в операциях сравнения используется операция сложения чисел, представленных в дополнительном коде. В качестве примера ниже показан процесс формирования значения деления $y_{15}^t$ шкалы $\mathbf{Y}_1^t$ и операция сравнения $y_{15}^t \; (><=) \; x$ для k = 4, q = 2, при этом x,z представляют мантиссы нормализованного числа в формате с плавающей запятой, т.е. $q^{-1} < x,z < 1$, $q^{-1} < y < q^1$.

| № | $y_{15}$ | $y_{15}^t = z\, y_{15}$ | $\Delta_{15} = y_{15}^t - x$ | $_{зн}\Delta_{15}$ |
|---|----------|--------------------------|------------------------------|---------------------|
| 15 | 1111 | $z_1\, z_2\, z_3\, z_4\, \ldots\, z_n$<br>$z_1\, z_2\, z_3\, z_4\, \ldots\, z_n$<br>$z_1\, z_2\, z_3\, z_4\, \ldots\, z_n$<br>$z_1\, z_2\, z_3\, z_4\, \ldots\, z_n$ | $z_1\, z_2\, z_3\, z_4\, \ldots\, z_n$<br>$z_1\, z_2\, z_3\, z_4\, \ldots\, z_n$<br>$z_1\, z_2\, z_3\, z_4\, \ldots\, z_n$<br>$z_1\, z_2\, z_3\, z_4\, \ldots\, z_n$<br>$1\; \bar{x}_1\, \bar{x}_2\, \ldots\, \bar{x}_n 1\; 1\ldots 1$<br>$\qquad\qquad\qquad\qquad\quad 1$<br>-------------------------------<br>$_{зн}\Delta\; \Delta_1 \Delta_2 \ldots\ldots\ldots\ldots \Delta_{n+k}$ | $\in 0,1$ |

Свертка пятирядного кода $\Delta_{15}$ в данном примере реализуется способом, описанным выше.

## 7. Распределенная арифметика в матричных процессорах

Совместим вычисление оператора суммы парных произведений **F** (5.1) с оператором группового суммирования **S** (2.1) в единой кодовой матрице **Q**:

$$\mathbf{Q} = \mathbf{F} + \mathbf{S} = \sum_{i=1}^{w}(A_i B_i + C_i + D_i + E_i + G_i + H_i + L_i) \qquad (7.1)$$

Положим формулу (7.1) в основу матричного арифметического процессора (МАП). Если на вход МАП поступают операнды в традиционном однорядном коде, то и величины A, B, C, D, E, G, H, L в формуле (7.1) рассматриваются как однорядные коды. Если на вход МАП поступают двурядные коды, например (C+D) или (E+H) или (A+H) и т.д., то и работа процессора рассматриваться в формате с двухрядными операндами.

Заметим, что при A=1 произведение AB =B (или при B=1 произведение AB=A) формула (6.1) будет представлять собой оператор группового суммирования семи однорядных кодов. На рис. 4 представлена обобщенная структура МАП, в основу работы которого положена формула (7.1). Данная структура рассчитана на операнды, имеющие разрядность не более 121 (n=< 121)

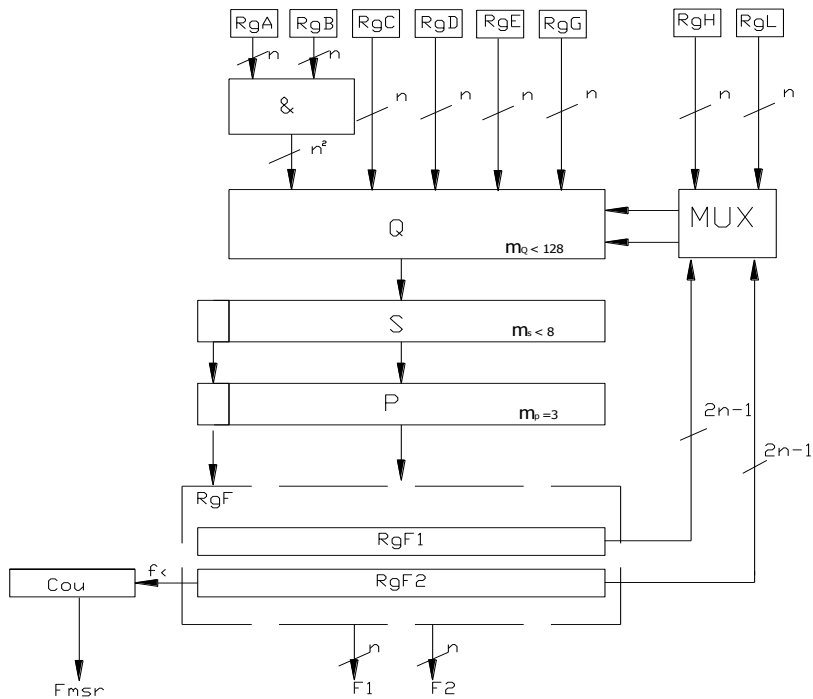

Fig. 4

В данной структуре матрица **&** представляет собой матрицу частичных произведений (матрица конъюнкций), матрица **Q** представляет собой матрицу сумматоров одноразрядных чисел, на входы которых сигналы поступают из матрицы **&** и входные операнды C, D, E, G, H, L. В режиме накопления вместо операндов H, L в матрицу **Q** с помощью мультиплексора MUX вводится двухрядный код (обратная связь) матрицы **F**. Свертка матрицы Q в двухрядный код F производится с помощью матриц **Q, S, P**, представляющих собой матрицы сумматоров одноразрядных чисел. Двухрядный регистр RgF (двухступенчатые регистры) запоминает результат свертки матрицы Q в формате двухрядного кода. Если все входные операнды представлены в прямом коде без знаковых разрядов (предполагается раздельная обработка положительных и отрицательных чисел), то проблема переполнения сетки (шириной 2n-1) решается путем суммирования всех разрядов, выходящих за пределы разрядной сетки, и счета переносов f (выход переноса из старшего сумматора одноразрядных чисел матрицы P) с помощью двоичного счетчика Cou (последовательного сумматора). Если входные операнды представлены в дополнительном коде (со знаковыми разрядами), то все матрицы формируются, исходя из ширины заданной разрядной сетки. В качестве примера на рис.5 показана форма матриц и распределение разрядов для 24-разрядных операндов A, B, C, D, находящихся в диапазоне (2 – 0 ). Время вычисления операнда F (6.1.) данном случае равно :

$$t_f = t_\& + t_q + t_s + t_p = 14 t_\& ,$$

где $t_\&$ - время формирования матрицы частичных произведений **&**, равное задержки в одном двухвходовом элементе AND;

$t_q = 6 t_\&$ - время свертки матрицы **Q** в матрицу **S** (см. таблицу 3.1.);

$t_s = 4 t_\&$ -время свертки матрицы **S** в матрицу **P**;

$t_p = 3 t_\&$ - время свертки матрицы **P**

Таким образом, оператор (6.1) на одной итерации для 24-х разрядных чисел вычисляется за $14 t_\&$.

Аппаратурные затраты для данной структуры ориентировочно составляют 12500 двухвходовых элементов AND.

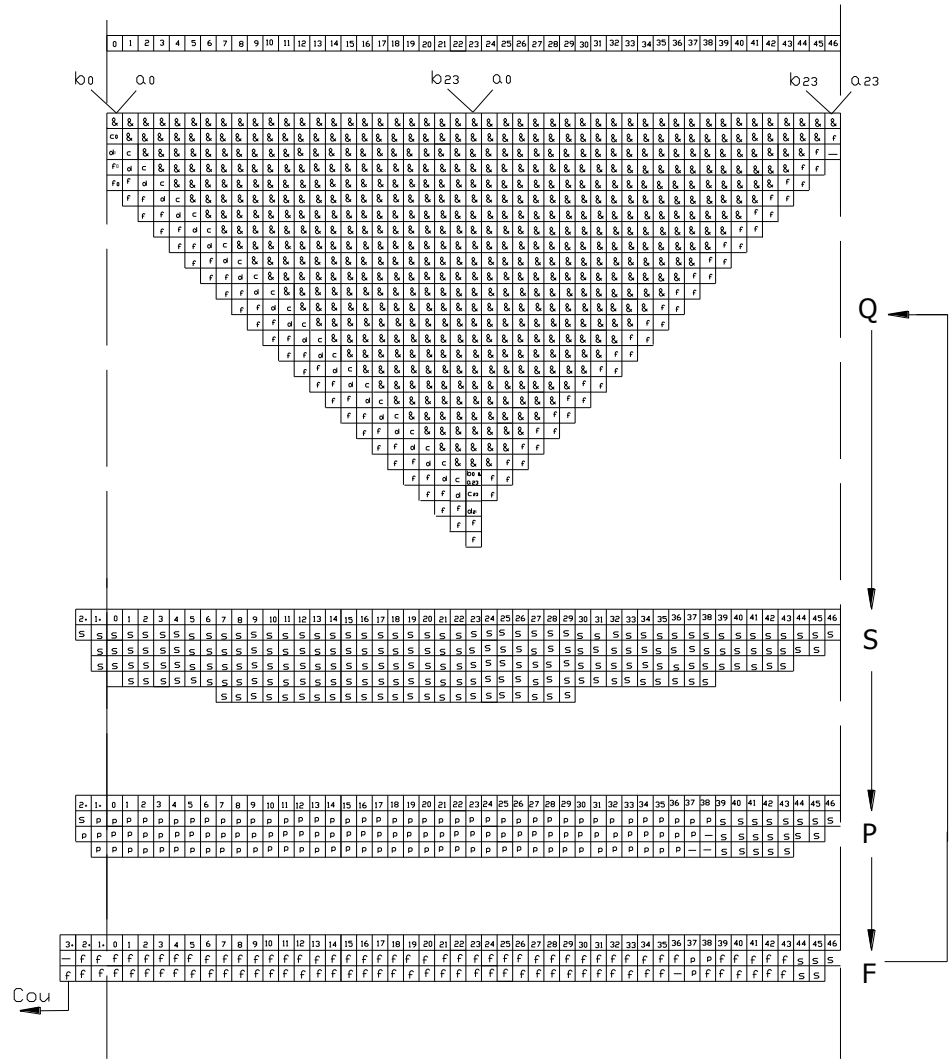

Fig. 5

.

## 8. Заключение

Вопросам повышения скорости выполнения арифметических операций на базе распределенной арифметики в последние десятилетия уделялось большое внимание. Особенно много работ в этой области посвящено структурам для вычисления разностных уравнений, используемых в цифровых фильтрах структурам для решения уравнений, используемых в спектральном анализе, структурам для решения матричных уравнений и т.д. [4,5,6]. Большинство этих работ базируется на использовании в матричных структурах стандартных микросхем. В данной работе рассмотрены матричные структуры применительно к конструированию самих микросхем, функционирующих на основе распределенной арифметики. Рассмотренные матричные структуры позволяют многократно повысить скорость выполнения арифметических операций по сравнению с классическими структурами, в которых основным элементов является параллельный сумматор. Реализация матричных структур вычислительных устройств связана с большими аппаратурными затратами, что тормозило их практическую реализацию. Однако большие успехи в интегральной технологии микросхем позволяют уже в настоящее время реализовать рассмотренные в данной статье идеи и по-новому подойти к архитектурам суперкомпьютеров. Следует так же учесть, что предлагаемые в данной работе матричные структуры имеют большую степень регулярности, что очень важно для интегральной технологии.

**Список литературы**

1. **Евдокимов В.Ф., Стасюк А.И., Щербаков В.И**. Матричные вычислительные устройства. «Наукова думка» , 1993. С. 1-152
2. **Нестеренко С.А., Паулин О.Н.** Построение обобщенной модели операции свертки многорядных кодов при цифровой обработке сигналов.//Технология и конструирование в электронной аппаратуре.-2008.-№1-С.20-26
3. **Храпченко В.М.** Методы ускорения арифметических операций, основанные на преобразовании многорядного кода//Вопросы радиоэлектроники. Сер. Электронная вычислительная техника.-1965.-Вып.8-С.121-144
4. **Mintzer,L** «FIR filters with the Xilinx FPGA» FPGA, 92 ACM/SIGDA Workshop on FPGAs hh.129-134
5. **Santoro M.R.** Design and clocking of VLSI multipliers/Stanford Universitey, Computer Systems Laboratory. Report Number: CSL-TR-89-397, October 1989.
6. The Role of Distributed Arithmetic in FPGA-based Signal Processing, Xilinx, Inc.1996

Robotic translation from Russian

*In this paper we consider the matrix structure of arithmetic processors based on distributed arithmetic in multi-row codes. Scope - development of supercomputers.*

*Key words: distributed arithmetic; multi-row codes; one-column adders; group sum; adder-accumulator; multipliers; vernier-conversion; arithmetic array processor.*

Ph.D.   **Shcherbakov V.I.,**  x162@rambler.ru
National Academy of Sciences of Ukraine
Pukhov Institute for Modelling in Energy Engineering

**Ultrafast a Distributed Arithmetic in multi-row codes**

**Introduction**

  Continually increasing demands on processing speed  leads to a need have a great computational resources. Modern supercomputers, representing an array of elementary computers, despite use the parallelization algorithms, they have big time costs  for transfer information between computers. Needed new ideas in the organization of computing processes, the creation of new architectures of supercomputers, based on advances in manufacturing technology of integrated circuits.  In this paper, we propose a method to reduce time costs associated directly with arithmetic processes.  The method is based on the algorithms in multi-row codes and using distributed arithmetic, which allows to parallelize computation process at the level of elementary operations with individual bits. We also offer the new architecture of the processors elements operating in multi-row codes focused on the integrated chip technology.

1. **The representation of the numbers in the multi-row codes**

  In modern computer technology used mainly positional codes, where each digit occupies a certain position in the code and has a certain significance**.** By summing the numbers, represented in floating-point format, you must align the order of the numbers, and therefore will be  lost bits mantissa  beyond the bit grid, and, consequently, the information is lost. Appear errors are associated with the loss of accuracy  in the consequence not a complete the intersection mantissas  on a number line.  Fixed-point format allows to summarize a big numbers with small numbers without loss of information, but it has range of representation of numbers, of limited bit grid .

By increasing the code word length is increased and the time of summation of numbers, so how should be carry-over down to most significant bit. Increased information of code can by the presentation operand in format matrix size m x n, where m - number of n- bit words. For example, a matrix of two n- bit words may be expressed as:

| $a_{1,n}$ | ... | $a_{1,i}$ | ... | $a_{1,1}$ |
|---|---|---|---|---|
| $a_{2,n}$ | ... | $a_{2,i}$ | ... | $a_{2,1}$ |

All codes contained in a particular column of this matrix have a certain ( equal ) weight. Subsequently, a matrix composed only of the codes will be **called matrix-codes or multi-row codes** . A remarkable property such the formats is that of the summation of matrices the carry is not spread along the entire bit-grid, and limited to a few bits. Suppose the numbers of A and B are defined in the 2-row codes:

|   | $a_{1,n}$ | $a_{1,n-1}$ | ... | $a_{1,i}$ | ... | $a_{2,1}$ | $a_{1,1}$ |
|---|---|---|---|---|---|---|---|
| A = | $a_{2,n}$ | $a_{2,n-1}$ | ... | $a_{2,i}$ | ... | $a_{2,2}$ | $a_{2,1}$ |
|   | $b_{1,n}$ | $b_{1,n-1}$ | ... | $b_{1,i}$ | ... | $b_{2,1}$ | $b_{1,1}$ |
| B = | $b_{2,n}$ | $b_{2,n-1}$ | ... | $b_{2,i}$ | ... | $b_{2,2}$ | $b_{2,1}$ |

Then the sum of these numbers will be a convolution ( 4→2) 4-row codes in 2-row code [3]. Can to note that the time summing will be in constant , which not independent of the number of digits (the number of bits of code). Denote the sum four one- bit numbers i -th column as $C_i$. If the code of the matrix A and B represented in binary codes, then Ci is represented by three bits with corresponding the weight: $C_i \to c_{i,3}\ c_{i,2}\ c_{i,1}$ .

Given that within column the position of bit not matter, then after the first phase summation of we obtain 3-row a code matrix of partial sums **C** (3-row code):

|   | - | - | $c_{n,1}$ | $c_{n-1,1}$ | ... | $c_{i,1}$ | ... | ... | $c_{2,1}$ | $c_{1,1}$ |
|---|---|---|---|---|---|---|---|---|---|---|
| C = | - | $c_{n,2}$ | $c_{n-1,2}$ | ... | $c_{i,2}$ | ... | ... | $c_{2,2}$ | $c_{1,2}$ | - |
|   | $c_{n,3}$ | $c_{n-1,3}$ | ... | $c_{i,3}$ | ... | ... | $c_{2,3}$ | $c_{1,3}$ | - | - |

As a result of convolution "3→2" we obtain the code 2-row matrix **S**, representing the sum of numbers A and B, in the same format as the input operands:

**S** =

| $s_{1,n+2}$ | $s_{1,n+1}$ | $s_{1,n}$ | ... | $s_{1,i}$ | ... | $s_{1,2} = c_{2,1}$ | $s_{1,1} = c_{1,1}$ |
|---|---|---|---|---|---|---|---|
| $s_{2,n+2} = c_{n,3}$ | $s_{2,n+1}$ | $s_{2,n}$ | ... | $s_{2,i}$ | ... | $s_{2,2} = c_{1,2}$ | $s_{2,1} = 0$ |

Thus, the sum of the two operands - matrices A and B size 2 x n gives matrix **S** of size 2 x (n + 2). At the same time, regardless of value n carry spreads only to the in three neighboring discharges. Time is determined by the temporary in one adder of four numbers one-column (in compressor 4→3) and one adder of three one-bit numbers (in compressor 3→2). The number A in the one-row code ( for example , fixed-point format) in a 2-row code will be as follows:

A =

| $a_{1,n}$ | ... | $a_{1,i}$ | ... | $a_{1,1}$ |
|---|---|---|---|---|
| 0 | ... | 0 | ... | 0 |

Size of matrix operands can be different and is limited only by the technological capabilities of manufacturing integrated circuits [1]. If the most high bits in of rows in the matrix m x n (high column) allocate for the sign bit, then such a code matrix b may representated additional code the numbers. Consequently, the subtraction operation in multi-row codes (as well as summing operation) are reduced to a convolution matrix, expressed as additional codes. When processing large data sets formation of additional codes requires large time costs. In addition, to additional code is required certain size of the bit grid, that when processing a large array of numbers spring up problem with overflow digit grid. For such problems we offer 4-row the matrix operands, in which the positive number locate in the two upper (or lower) rows, and negative in the two lower (or upper) rows:

A =

| | | $a_{1,n}$ | $a_{1,n-1}$ | ... | $a_{1,i}$ | ... | $a_{2,1}$ | $a_{1,1}$ |
|---|---|---|---|---|---|---|---|---|
| $A_+$ = | | $a_{2,n}$ | $a_{2,n-1}$ | ... | $a_{2,i}$ | ... | $a_{2,2}$ | $a_{2,1}$ |
| | | $a_{1,n}$ | $a_{1,n-1}$ | ... | $a_{1,i}$ | ... | $a_{2,1}$ | $a_{1,1}$ |
| $A_{(-)}$ = | | $a_{2,n}$ | $a_{2,n-1}$ | ... | $a_{2,i}$ | ... | $a_{2,2}$ | $a_{2,1}$ |

In this format positive numbers are added to the positive numbers, and negative to negative. The result of the summation is also presented in the same format.

$$S = A + B = \begin{array}{|c|} \hline A_{(+)} \\ \hline A_{(-)} \\ \hline \end{array} + \begin{array}{|c|} \hline B_{(+)} \\ \hline B_{(-)} \\ \hline \end{array} = \begin{array}{|c|} \hline S_{(+)} = A_{(+)} + B_{(+)} \\ \hline S_{(-)} = A_{(-)} + B_{(-)} \\ \hline \end{array}$$

When the subtraction, then in subtrahend must invert matrices (top to bottom source code) and then to perform the addition:

$$S = A - B = \begin{array}{|c|} \hline A_{(+)} \\ \hline A_{(-)} \\ \hline \end{array} + \begin{array}{|c|} \hline B_{(-)} \\ \hline B_{(+)} \\ \hline \end{array} = \begin{array}{|c|} \hline S_{(+)} = A_{(+)} + B_{(-)} \\ \hline S_{(-)} = A_{(-)} + B_{(+)} \\ \hline \end{array}$$

It is obvious that in the subtraction in this case does not required additional-tional codes. All codes $A_{(+)}$, $A_{(-)}$, $B_{(+)}$, $B_{(-)}$, $S_{(+)}$, $S_{(-)}$ in this format no sign bits and their formats this is direct codes. The weight of the bits in junior column of the matrix is assigned a special code (companion matrix) whose value is a multiple byte. This allows the position of the point that separates the integer part from the fractional part, define virtual by meaningfully weight discharges junior column. This format can be termed as «format with a virtual position of the point». For example, if the weight of discharges in the junior column is set to $2^{-128}$, the point will be located at 127 bits (columns) from junior column. In this representation the operands operation accumulates can be made without loss of accuracy if the weight bits in the junior column equal to the minimum possible value of the input number. Data bus for input and output data in this format should be a 4-bus line. All arithmetic operations (addition, subtraction, multiplication, division), ultimately boil down to one type of operation - the operations of addition, and in the proposed format to the convolution multi-row code. Consider the calculation of the operator group summation $S = \Sigma A_i$ and operator of a sum of product pair of $F = \Sigma A_i B_i$ (i = 1 ... m). These operators form the basis for calculating the difference equations, equations of the "butterfly" in the spectral analysis, matrix equations, the equations for the digital filters, etc. The speed of these operators largely determine the speed of the entire computational process in the processing of the data stream.

## 2. Computation the operator of the group summation

Operator group summation

$$S = \sum_{i=0}^{m-1} A_i \quad , \quad A_i = \sum_{j=0}^{n-1} a_{i,j} g^j \tag{2.1}$$

may be expressed in matrix form as:

$$S = A\,Q$$

or

$$S = \begin{array}{|c|c|c|c|c|} \hline a_{0,n-1} & \ldots & a_{0,j} & \ldots & a_{0,0} \\ \hline \ldots & \ldots & \ldots & \ldots & \ldots \\ \hline a_{i,n-1} & \ldots & a_{i,j} & \ldots & a_{i,0} \\ \hline \ldots & \ldots & \ldots & \ldots & \ldots \\ \hline a_{m-1,n-1} & \ldots & a_{m-1,j} & \ldots & a_{m-1,0} \\ \hline \end{array} \quad \begin{array}{|c|} \hline g^{n-1} \\ \hline \ldots \\ \hline g^j \\ \hline \ldots \\ \hline g^0 \\ \hline \end{array}$$

Matrix **A** which formed from m one-row code, called the code-matrix or m-row code. Imagine operator (2.1) in redundant number system as:

$$S = \sum_{i=0}^{m-1} \sum_{j=0}^{n-1} a_{i,j} g^j = \sum_{j=0}^{n-1} S_j\, g^j \tag{2.2}$$

where $S_j = \sum_{j=0}^{m-1} a_{i,j}$ - the sum of the one-bit-digit numbers in j-th column of matrix **A**

Since m-row the code S need reduce up 2-row code, then summation process will take place in several stages (layers). The sum of S, as well as matrix-codes corresponding to the stage of the convolution will be denoted by the index in parentheses. Thus, the original code matrix **A** in expression (2.1) will be designated **A** $_{(1)}$ and the sum of the one-bit numbers in the j-th column of this matrix will have $S_{(1),j}$. Similarly, for the dimensions of the matrix **A** $_{(1)}$ : m → m $_{(1)}$ ; n → n $_{(1)}$ …

Representation the S (1), j in the number system with radix q:

$$S_{(1),j} = \sum_{h=0}^{m_{(2)}-1} a_{i,j}\, q^h \tag{2.3}$$

where $m_{(2)} = ]\log_q (1+S_{(1),j}{}^{max})[ \;=\; ]\log_q (1+m_{(1)}q - m_{(1)})[$ , (2.4)

$S_{(1),j}{}^{max} = m_{(1)}(q-1)$ - the maximum value of the amount of one-bit numbers in the j-th column; ] ... [ - notation rounding up to the nearest integer.

We write the expression (2.2) considering (2.3) in matrix form

$$S = A_{(2)} Q \; ,$$

wherein the code matrix $A_{(2)}$ has the following form:

$A_{(2)} = $

|  |  |  |  |  | $a_{0,n-1}$ | ... | $a_{0,j}$ | ... | $a_{0,1}$ | $a_{0,0}$ |
|---|---|---|---|---|---|---|---|---|---|---|
|  |  |  |  | $a_{1,n-1}$ | ... | $a_{1,j}$ | ... | $a_{1,1}$ | $a_{1,0}$ |  |
|  |  |  | ... | ... | ... | ... | ... | ... |  |  |
|  |  | $a_{i,n-1}$ | ... | $a_{i,j}$ | ... | $a_{i,1}$ | $a_{i,0}$ |  |  |  |
|  | ... | ... | ... | ... | ... | ... |  |  |  |  |
| $a_{m-1,n-1}$ | ... | $a_{m-1,j}$ | ... | $a_{m-1,1}$ | $a_{m-1,0}$ |  |  |  |  |  |

The matrix $A_{(2)}$: $n = m_{(2)}$ ; $m = n_{(1)}$

As matrix $A_{(2)}$ as well as matrix $A_{(1)}$ reflect the sum of S, then the location of the bits $a_{i,j}$ within the j-th column is irrelevant. Therefore, shifting the corresponding bits within their own columns in a matrix form the matrix $A_{(2)}$ in the form of a trapezoid:

$A_{(2)} = $

| $a_{m-1,n-1}$ | ... | $a_{i,n-1}$ | ... | $a_{1,n-1}$ | $a_{0,n-1}$ | ... | $a_{0,j}$ | ... | $a_{0,1}$ | $a_{0,0}$ |
|---|---|---|---|---|---|---|---|---|---|---|
|  | ... | ... | ... | ... | ... | ... | ... | ... | ... |  |
|  |  | ... | ... | ... | ... | ... | ... | ... |  |  |
|  |  |  | ... | ... | ... | ... | ... |  |  |  |
|  |  |  | $a_{m-1,j}$ | ... | $a_{m-1,0}$ |  |  |  |  |  |

In this matrix: $m = m_{(2)}$; $n = n_{(1)} + m_{(2)} - 1$.

Code length in the (m-1) th row trapezoidal the matrix $A_{(2)}$ is minimum, and is equal to $n_{min} = n_{(1)} - m_{(1)} + 1$. The matrix $A_{(2)}$ is $m_{(2)}$ -row code, which correlate with $m_{(1)}$ - row code of logarithmic the dependence (2.4). A feature of the code matrix $A_{(2)}$ is that it also reflects of matrix partial products of two one-row codes.

For example, product $F = BC$, where $B = b_{n-1} q^{n-1} + ... b_i g^i + ... b_0 q^0$, $C = c_{n-1} q^{n-1} + ... c_i q^i + ... c_0 q^0$, (i-positive integers) can be represented in the form of a code F :

|  |  |  |  | $b_0 \& c_{n-1}$ | ... | $b_0 \& c_i$ | ... | $b_0 \& c_0$ |
|---|---|---|---|---|---|---|---|---|
|  |  |  | ... | ... | ... | ... | ... |  |
| $\mathbf{F}_{(2)} =$ |  |  | ... | ... | ... | ... | ... |  |
|  |  | ... | ... | ... | ... | ... |  |  |
|  |  | ... | ... | ... | ... |  |  |  |
|  | $b_{n-1} \& c_{n-1}$ | ... | $b_{n-1} \& c_i$ | ... | $b_{n-1} \& c_0$ |  |  |  |

Therefore, all subsequent stages of the convolution multi-row code will reflect the process of multiplying two one-row codes.

Each subsequent matrix-code $\mathbf{A}_{i+1}$ formed in the process of convolution from the matrix-code $\mathbf{A}_i$ using addition operation one-bit numbers (2.3). Wherein the number of rows followed by a matrix linked with previous the matrix logarithmic dependence (2.4). Table 2.1 shows the numbers row $m_i$ at each stage, by convolution of the original m-row the code (for the binary system).

Таблица 2.1

| $m_{(1)}$ - the number of rows | $m_{(2)}$ | $m_{(3)}$ | $m_{(4)}$ | L - number of stage |
|---|---|---|---|---|
| $2^2 - 1$ | 2 | - | - | 1 |
| $2^2 ... (2^3 - 1)$ | 3 | 2 | - | 2 |
| $2^3 ... (2^4 - 1)$ | 4 | 3 | 2 | 3 |
| $2^4 ... (2^5 - 1)$ | 5 | 3 | 2 | 3 |
| $2^5 ... (2^6 - 1)$ | 6 | 3 | 2 | 3 |
| $2^6 ... (2^7 - 1)$ | 7 | 3 | 2 | 3 |

Analysis of the data in Table 2.1 indicates that is a step function $L = f(m_i)$, where L- number of stages in convolution the original matrix-code in 2-row code. Thus the value L in a certain range does not depend on the number of addends. For example, 64 or 127 addends can compressing in of 2-row code in 3 stages This fact (in constructing the apparatus) indicates the possibility to speed up the process computing by means of selection optimal sizing of source matrix-code. Hardware costs for multi-row convolution codes are mainly determined by the number (and the number of inputs) one-column adders (further – OCA).

Calculate these costs can based on the following relationships:

for matrix $\mathbf{A}_{(2)}$ :

$n_{(2),i} = n_{(2)}{}^{max} - \Delta n_{(2),i}$ - code length i-row

$\Delta n_{(2),i} = 0$, if i=0;

$\Delta n_{(2),i} = 2^i$, if $i < m_{(2)}$, i= 1, …($m_{(2)} - 1$) − number row

$N_{(2),v} = 2$, if $v < m_{(2)}$, - the number of columns having a length v

$N_{(2),v} = n_{(2)}{}^{min}$, if $v = m_{(2)}$, v = (1, … $m_{(2)}$) − column length (number of bits)

$n_{(2)}{}^{max} = n_{(1)} + n_{(2)} - 1$, $n_{(2)}{}^{min} = n_{(1)} - m_{(2)} + 1$

for matrix $A_{(3)}$:

$n_{(3),i} = n_{(3)}{}^{max} - \Delta n_{(3),i}$ - code length i-row

$\Delta n_{(3),i} = 0$, if i=0;

$\Delta n_{(3),i} = 2^i$, if $i < m_{(3)}$, i= 1, …($m_{(3)} - 1$) − number row

$N_{(3),v} = 2$, if $v < m_{(3)}$, - the number of columns having a length v

$N_{(3),v} = n_{(3)}{}^{min}$, if $v = m_{(3)}$, v = (1, … $m_{(3)}$) − column length (number of bits)

$n_{(3)}{}^{max} = n_{(2)}{}^{max}$

$n_{(3)}{}^{min} = n_{(3)}{}^{max} - q^{m(3)} + 1$

The time costs to compute the operator group summation is mainly determined by the speed summing in OCA and the number of stages of the convolution. In turn, the performance of OCA is determined by the number addends. Obviously, the greatest number of addends ($m_{(1)}$) in the same column comprises a matrix $A_{(1)}$.

The compression time of the original matrix size n x $m_1$ in 2-row the code determined by the number of stages convolution L, and the time $t_m$ summing in the OCA in each matrix, which, depends on the number of addends m.

For example, if the number of addends ( for any n) in the original matrix does not exceed 63 ($m_1 = 63$), then time the matrix-code compression in 2-row code using tabular - OCA will be:

$$t_{mxn} = t_{m1} + t_{m2} + t_{m3} + t_{cc} = 13t_\&$$ (2.5)

$$t_{m1} = 6t_\&, \ t_{m2} = 3t_\&, \ t_{m3} = t_\&, \ t_{cc} = 3t_\&$$

where $t_\&$ - a delay in one component AND,

$t_{cc}$ - a delay in encoder adder

$m_1=63, m_2 =6, m_3 = 3, m_4 =2$

### 3. The tabular one-column adder (tabular - OCA)

From the analysis of matrices-code $A_{(1)}, A_{(2)}, A_{(3)}$, that the minimum number of inputs in the OCA this is 3, and the maximum is equal to the number of addends in the original matrix ($m_{(1)}$). If $m_{(1)}$ is large, then implement the tabular-OCA (within the allowable hardware costs) can be in the form of a tree structure [1]. Below we will consider only tabular adders , because their structure more fully meets the requirements of IC. Consider an elementary tabular - OCA having three inputs and two outputs .

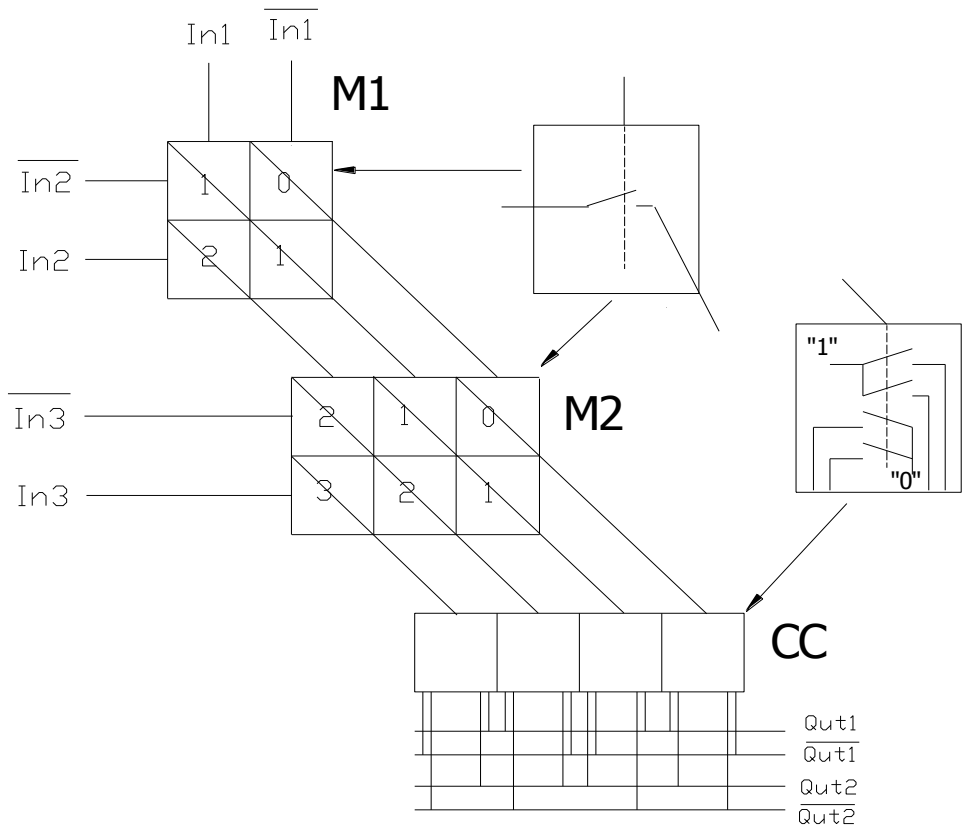

Fig.1

Fig. 1 shows the structure of the tabular - OCA 3→ 2. Componentry in table M1, M2 this two-input elements AND.. Elements encoder - CC can be performed on the transistors of type "switch". Note that the table elements M2, M3 may be in the form of the switch. A first input (e.g., gate) fed an enable signal and a second input (e.g., a drain or source) for enable information signal. Timed delay in the tabular - OCA is equal to $t_{3-2} = 2 t_\& + t_k$, where $t_\&$, $t_k$ is timed delay in the in two-input AND element and one of the switch. If the entire OCA is built on switches then: : $t_{3-2} = 3t_k$. For all this hardware costs account for 26 switch.

Consider tabular OCA, having the m inputs . Such OCA will have a tree structure. To simplify the explanation, we will consider only square tables that make up the tree structure OCA, which shown in Fig. 2 . The tree structure tabular OCA for the 64 numbers contains 6 types of tables: M1 ... M6. In the general case (of the square table) number of types of tables N related to the number of addends m by equality: $m = 2^N$. For all this table have the following dimensions: M1 = 2 x 2; M2 = 3 x 3 ; M3 = 5 x 5 ; M4 =9x9; M5 = 17 x 17; M6 = 33 x 33. These dimensions (for square tables) do not depend on the number of addends. Number of tables $V_i$, where i is the number of types table, related to the number (m) of addends (divisible by $2^i$) the following by equation: $V_i = m / 2^i$.
For example, for m = 64, we have: M1 = 32; M2 = 16; M3 = 8; M4 = 4; M5 = 2; M6 = 1. Since the table, (for any dimension) has temporary delay equal $t_\&$ where $t_\&$ be temporary delay in one the component AND ( with two-input) , that total time of summation of 64 OCA will be equal to $t_s = 7 t_\&$, including the delay in encoder CC $t_{cc} = t_\&$. In the range the inputs (3-128) time addition m - one-bit numbers can be determined from Table 3.1

Table 3.1

| Number of inputs | 3-4 | 5-7 | 8-16 | 17-32 | 33-64 | 65-128 |
|---|---|---|---|---|---|---|
| Time summation | $2t_\& + t_{cc}$ | $3t_\& + t_{cc}$ | $4t_\& + t_{cc}$ | $5t_\& + t_{cc}$ | $6t_\& + t_{cc}$ | $7t_\& + t_{cc}$ |

The table shows that the time of summation is expressed by a step function that must be considered when designing devices for convolution of multi-row codes. Hardware costs tabular OCA for range of inputs (3-32) are summarized in Table 3.2.

Table 3.2

| M, $\Sigma_\&$ | The number of inputs of the adder (m) | | | | | | | | | | | | | | |
|---|---|---|---|---|---|---|---|---|---|---|---|---|---|---|---|
| | 3 | 4 | 6 | 8 | 10 | 12 | 14 | 16 | 18 | 20 | 22 | 24 | 26 | 28 | 30 | 32 |
| M1 | 1 | 2 | 3 | 4 | 5 | 6 | 7 | 8 | 9 | 10 | 11 | 12 | 13 | 14 | 15 | 16 |
| M2 | 1 | 1 | 1 | 2 | 2 | 3 | 3 | 4 | 4 | 5 | 5 | 6 | 6 | 7 | 7 | 8 |
| M3 | - | - | - | - | 1 | 1 | 1 | 2 | 2 | 2 | 3 | 3 | 3 | 3 | 4 | 4 |
| M4 | - | - | - | - | 1 | 1 | 1 | 1 | 1 | 1 | 1 | 1 | 2 | 2 | 2 | 2 |
| M5 | - | - | - | - | - | - | - | - | 1 | 1 | 1 | 1 | 1 | 1 | 1 | 1 |
| $\Sigma_\&$ | 10 | 17 | 36 | 59 | 90 | 121 | 158 | 199 | 254 | 301 | 354 | 411 | 476 | 541 | 582 | 687 |

This table shows the number of matrices type M1-M5 and the total number of two-input elements AND ($\Sigma_\&$), incoming in of encoders CC, for a given number of inputs (m) OCA. For 64 numbers ( one-bit) of quantity of two-input elements AND in the encoder CC will be equal to $\Sigma_\& = 2463$.
Convolution of multi-row codes based on standard chips is described in [2, ].

## 4. Group summing of the accumulation

We write group summing of the accumulation in the following form:

$$S = \sum_{i=1}^{w} S_i, \qquad S_i = \mathbf{A}_i + S_{i-1} \qquad (4.1)$$

where $S_{i-1}$ - the current value of the sum; $\mathbf{A}_i$ - matrix - operand input.

In this case, at each iteration is input: operand - matrix $\mathbf{A}_i$, having the size of the array (n x m), and the result convolution of the preceding multi-row code $S_{i-1}$. Given that (for convolution multi-row code) 2-row code is always preceded by a 3-row коде, then in the expression (4.1) $S_{i-1}$ can represent 3-row code. This reduces the number of steps of the convolution at one step and accelerate the process of convolution. On (w-1) - iteration must compressing 3-row code $S_{w-1}$ in 2-row code. Consider the most simple architecture of the adder- accumulator (Fig.3), in which the operand $A_i$ is represented in 1-row code, and the current sum $S_{i-1}$ is 2-row code.

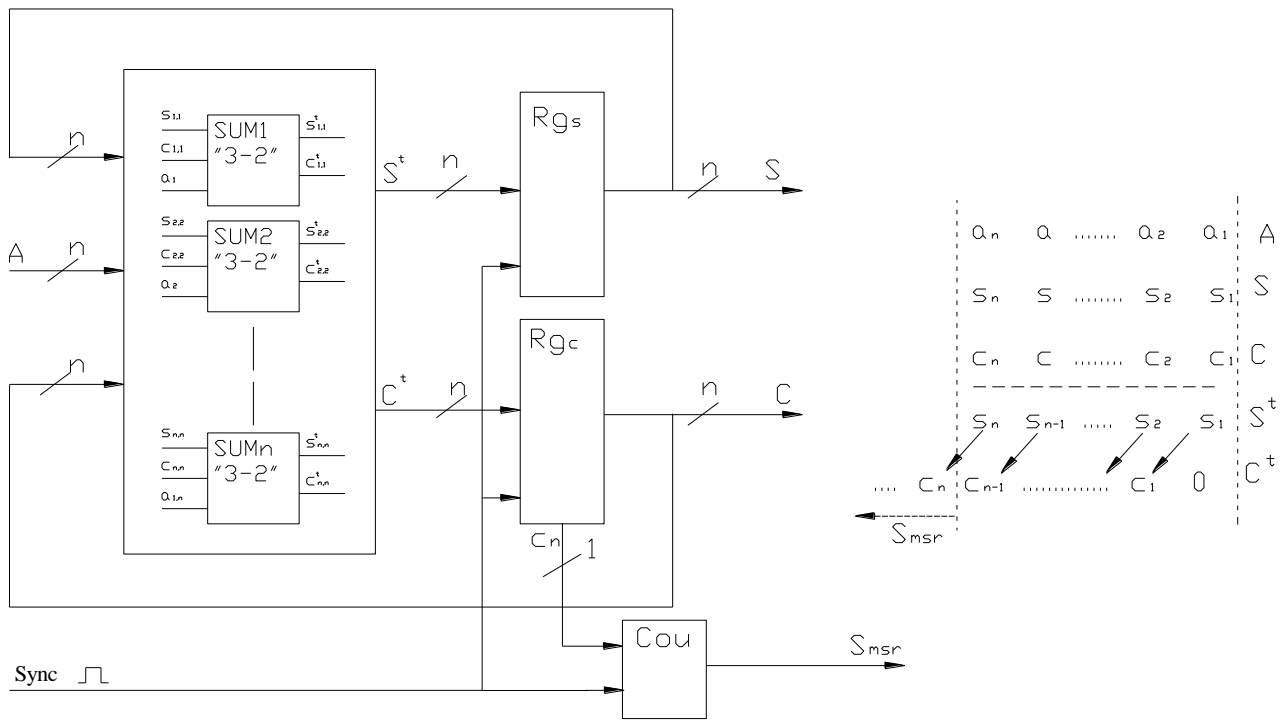

Fig. 3

Supposed to separate summation of positive and negative numbers. All numbers are represented in direct codes without sign bits. In this case the basis of the adder-accumulator will be ruler from n - adders OCA such as $S_{3 \rightarrow 2}$, the role of which can carry a full - adder. The structure of the adder also includes two registers Rgs, Rgc (triggers type "master-slave"). Register Rgs receives from the adders bits of sum and register Rgc - bits of carry . Location of the bits in the summation is shown in the same figure in the form of a callout. Time of summation n-bit numbers Ai is determined by the time of summation in of one adder $S_{3 \to 2}$ and the time of recording information in the register.    Since the carry-over of the summation in this structure carry-over only on to one bit, then all the input and output of codes can be placed in the n-bit grid. Overflow in summation determined by as the carry bit $c_n$ from n-th (MSB) adder $S_{3 \to 2}$. This carry arrives at the input to the binary counter impulses Cou with parallel binary output. In such a structure, an adder-accumulator can have arbitrarily large range of the output code and add up the numbers without losing precision.   Moreover, if, for example, is required to accumulate the of low-bit-numbers $A_i$ with a very large number of iterations (cycles of accumulation w), then such a structure will have a minimum hardware costs because the lines OCA and registers also will have a small length the code.    In this case, code length register Rgs must be increased by one bit. Output bit $s_{n+1}$ ( from Rgs) must be addition modulo-two with the output bit $c_{n+1}$ , (from Rgc) and summation result enter at input of counter Cou .

### 5.   The arithmetic sum of pair products

Before considering the process of calculating the operator of sum of pair products

$$F = \sum_{i=1}^{m} A_i B_i \qquad (5.1)$$

consider the product P  of numbers A, B represantation in the additional code :

$$P = AB,$$

where  $A = -a_{(0)}q^0 + \sum_{i=1}^{n} a_i q^{-i}, \quad B = -b_{(0)}q^0 + \sum_{j=1}^{n} b_j q^{-j},$

$a_{(0)}$, $b_{(0)}$ - the sign bit, and numbers of A, B can be normalized or denormalized.

We write the product of these numbers as :

$$P = AB = D + E - F - C, \qquad (5.2)$$

where  $D = a_{(0)}b_{(0)}q^{(0)}$ ;  $E = \sum_{i=1}^{n}\sum_{j=1}^{n} a_i b_j q^{-(i+j)}$ ;  $F = \sum_{j=1}^{n} a_0 b_j q^{-j}$ ;  $C = \sum_{i=1}^{n} b_0 a_i q^{-i}$

The addend E representative of matrix of partial productsю and does not depend on sign bits while addends D, F, C - corrective the numbers for computation of product in the additional code. If $a_0 = 0$ we have $D = F = 0$, and if $b_0 = 0$, we obtain $D = C = 0$.

Replace the subtraction in expression (5.2) on the summation, for this represent the number of F and C in the additional code:

$$P = AB = D + E + F_{ac} + C_{ac} \qquad (5.3)$$

Code matrix P, formed in accordance with the expression (5.3), has the following form:

| $2^0$ | $2^{-1}$ | $2^{-2}$ | $2^{-3}$ | ... | ... | ... | ... | ... | ... | $2^{-2n}$ |
|---|---|---|---|---|---|---|---|---|---|---|
| $b_0a_0$ | **$b_1a_0$** | $b_1a_1$ | $b_2a_1$ | $b_3a_1$ | ... | ... | ... | ... | ... | $b_na_n$ |
| | **$b_0a_1$** | $b_2a_0$ | $b_1a_2$ | $b_2a_2$ | $b_2a_3$ | ... | ... | ... | ... | |
| | | **$b_0a_2$** | $b_3a_0$ | ... | ... | ... | ... | ... | | |
| | | | **$b_0a_3$** | ... | ... | ... | ... | | | |
| | | | | ... | **$b_n\,a_0$** | ... | | | | |
| | | | | 1 | **$b_0a_n$** | | | | | |

Bold type in the matrix shows the inverted AND. Matrix P, this is the matrix of partial products corresponding to the matrix A (2), discussed in Section 2. Therefore the process of convolution the matrix-row code and estimation of temporary costs and hardware costs corresponds to convolution the matrix A (2) which discussed above. The number of stages L convolution matrix P can be determined from Table 2.1. Note that if the operands A and B representated in straight codes without sign bits, then in the above matrix of partial products is necessary to exclude unity in the last line and the inverted bits replace on are not inverted bits, wherein of sign bits $a_0$, $b_0$ be the information bits with a weight of $2^0$. Time multiplying will be equal to the delay of formation of the matrix P, equal delay $t_\&$ in single element AND, and time convolution matrix P in the 2-row code. For example, time multiplying two 63-bit codes (look 2.5) will be equal to 14 $t_\&$.

The represent the operator of a sum of product pairt (5.1) as a recurrent formula:

$$F_i = F_{i-1} + A_iB_i = F_{i-1} + \mathbf{P}_i$$

Due to the fact that a function $L = f(m_i)$ this a step function, then operand $F_{i-1}$ can be included in the code matrix **P**. This eliminates the addition operation in the given expression. In order not to increase the number of stages of the convolution ( and time convolution respectively) of the matrix **P**, it is necessary to determine the stage in which $F_{i-1}$ be included.   Note that the code matrix **P** can be convolution in a 2-row code by   the multiplicators tabular (based narrow  width of memory) [1] or on standard chip [2,3,5].

## 6.  Operation division

The operation of division can be written in the following form:

$$y = x / z,$$

or   $y^t = x$,   where   $y^t = yz$ \hspace{4em} (6.1)

From the expression (6.1) it follows that in   operation division the involved operations  multiplication and comparison. In  modern algorithms implementing the division operation, applies a method of calculating "digit by digit" in which  on each iteration is determined one bit private.  Consider  the accelerated method of division, wherein on each iteration are calculated k bits of quotient the relation. The value of k is limited only by the technological capabilities of the integrated technology.   The essence of the method consists in of vector comparison operation x with a digital scale x $Y_1^t$ on the first iteration and calculating the difference ( discrepancy ), after which in a comparison operation on the second iteration  used discrepancy, etc.  Let us consider this process in more detail.   We introduce the concept of linear digital scale. Under   digital scale of  **Y** we understand the a linear scale having $q^n$ - digital point (including the zero point), where q-radix.   All weight the values points  presented in n-bit the codes.  The value of each of the two neighboring point at such a scale different from each other on one "quantum", ie the on value of the least significant bit of the code. Thus, on such a scales would be present all the possible values of n-bit number, and hence, the determination of the unknown n-bit  number *x* by comparing it with the scale, can be made in a single iteration.  Since this scale should be stored in memory or formed during the computation, then for a large-bit number the  hardware costs  will be unacceptably high.  Instead of one "long" digital scale **Y**, having $q^n$ points: $\mathbf{Y} \rightarrow y_1\ y_2 \ldots yi \ldots y_n$, in comparison operations can be using a set of digital scales in which the "length" of each subsequent scale overlaps the interval between two neighbors points of the previous scale.

$$\mathbf{Y} \rightarrow (y_0\ y_1\ \ldots\ y_k)_1\ (y_{k+1}\ y_{k+2}\ \ldots\ y_{k+2k})_2\ \ldots\ (y_{m-k}\ \ldots\ y_n)_m \rightarrow \mathbf{Y}_1\ \mathbf{Y}_2 \ldots \mathbf{Y}_j \ldots \mathbf{Y}_m$$

Thus, if the generated first scale $Y_1$, then all other vernier scales can generated by shifting codes points scale on the k-bits to the right at the beginning of each subsequent iteration, i.e.

$$Y = \sum q^{-k(j-1)} Y_1, \text{ where } j = 1, 2....m \qquad (6.2)$$

From the expression (6.1) it follows that is involved in a comparison operation is not digital scale $Y$, but its modification $Y^t$, which is the product of a scalar quantity z on the vector $Y$. I.e. all the points scale $Y$ to be multiplied on the n-bit code the number z. And the values (weight) of points of scales $Y$ and $Y^t$ a related of tabular (through a number of points). A comparison of the scalar variable $x$ with the vector $Y_1^t$ (6.1) on the first iteration we obtain the vector difference $\Delta$. By analyzing the sign bit $_s\Delta$, we obtain a unitary code that uniquely identifies the number point on the scale $Y_1^t$ and, accordingly, on the scale of $Y$. On the output of the decoder unitary code will form the first k high bits of quotient relation:

$$x \ (><=) \ Y_1^t \to \Delta_1 \to \Delta_{1,h} \to y^t_{1,h} \to y_{1,h} \to (y_0 \ y_1 \ .... \ y_k)_1$$

where h is the index number of the nearest point (on shortage), weight which approximately equal to $x$.

Simultaneously with decoding of a unitary code is determined and the value of discrepancy (on shortage):

$$\Delta_1 = y^t_{1,h} - x$$

On the second iteration with the scale $Y_2^t$ compared is discrepancy $\Delta_1$. Instead of the formation of scale $Y_2^t$ by means shifting its codes points to the right by k bits (6.2) can be shifted codes $\Delta_1$ on k bits to the left. The mathematical model of the process of division on discrepancies can representation in the following form:

$$\Delta_{j-1,h} \ (><=) \ Y_1^t \to \Delta_j \to \Delta_{j,h} \to y^t_{j,h} \to y_{j,h} \to (y_0 \ y_1 \ .... \ y_k)_j$$

$$\Delta_{j,h} = y^t_{j-1,h} - \Delta_{j-1} \quad \text{for } j = 1...m$$

$$\Delta_{j+1,h} = q^k \Delta_{j,h} \quad \text{for } j > 1$$

$$\Delta_{j-1,h} = x \quad \text{for } j = 1$$

Multiplication and comparison can be performed simultaneously using the convolution of multiple-row code. Wherein, instead subtraction in comparison operations used the operation of addition of numbers represented in the additional code. As an example below shows the process of forming values point $y_{15}^t$ scale $Y_1^t$ and comparison operation $y_{15}^t$ (> <=) x for k = 4, q = 2, wherein x, z is are normalized mantissa of floating-point, ie, . $q^{-1}$ <x, z <1, $q^{-1}$ <y <$q^1$.

| № | $y_{15}$ | $y_{15}^t = z\, y_{15}$ | $\Delta_{15} = y_{15}^t - x$ | $_s\Delta_{15}$ |
|---|---|---|---|---|
| 15 | 1111 | $z_1\ z_2\ z_3\ z_4\ ...\ z_n$<br>$z_1\ z_2\ z_3\ z_4\ ...\ z_n$<br>$z_1\ z_2\ z_3\ z_4\ ...\ z_n$<br>$z_1\ z_2\ z_3\ z_4\ ...\ z_n$ | $z_1\ z_2\ z_3\ z_4\ ...\ z_n$<br>$z_1\ z_2\ z_3\ z_4\ ...\ z_n$<br>$z_1\ z_2\ z_3\ z_4\ ...\ z_n$<br>$z_1\ z_2\ z_3\ z_4\ ...\ z_n$<br>$1\ \bar{x}_1\ \bar{x}_2\ ......\ \bar{x}_n 1\ 1...1$<br>$\qquad\qquad\qquad\qquad\qquad 1$<br>------------------------------<br>$_s\Delta\ \Delta_1 \Delta_2\ .................\ \Delta_{n+k}$ | € 0,1 |

Convolution 5-row code $\Delta_{15}$ in this example is implemented in the manner described above.

### 7. Distributed arithmetic in matrix processors

Compatible computing the sum of pair products F (5.1) with group summation **S** (2.1) in a single code matrix **Q**:

$$\mathbf{Q} = \mathbf{F} + \mathbf{S} = \sum_{i=1}^{w}(A_i B_i + C_i + D_i + E_i + G_i + H_i + L_i) \qquad (7.1)$$

Let the formula (7.1) be as the basis of the arithmetic matrix -processor (MAP). When the MAP receives input operands in a traditional one-row code, then and the values A, B, C, D, E, G, H, L in formula (7.1) are considered as one-row codes.

If on input receives MAP 2-row codes, such as (C + D) or (E + H) or (A + H) etc., then and functioning the processor be in format 2-row operands. Note that if A = 1 product AB = B (or B = 1 product AB = A), formula (6.1) will be the operator of the group summation of the 7-row codes. Fig. 4 is a generalized structure of the MAP, based functioning on a formula (7.1). This structure is designed for operands having a bit width no more than 121 (n = <121)

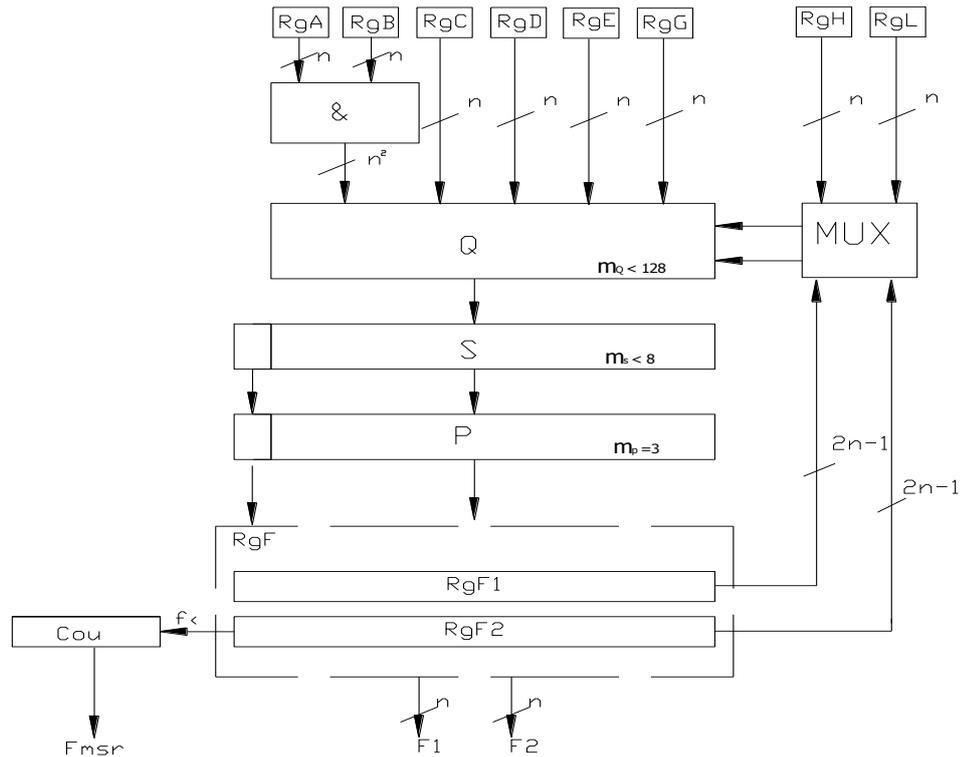

Fig. 4

In this structure, the matrix **&** this is a matrix of partial products (matrix conjunctions), the matrix **Q** is a matrix of adders one-row numbers on input which signals are received from the matrix **&** and the input operands C, D, E, G, H, L. In accumulation mode instead of operands H, L in the matrix **Q** via the multiplexer MUX is entered 2-row code (feedback) of the matrix F. Convolution matrix **Q** in 2-row code F produced with the help matrix **Q**, **S**, **P**, representing an arrays OCA. Registers RgF (two-stage registers) stores the result of the convolution matrix **Q** in format 2-row code. If all input operands are shown in direct code without sign bits (assuming separate processing of positive and negative numbers), then problem of overflow word length (length of 2n-1) eliminated with the by summation of bit grid overflow and summation carry f (f- carry output of the older OCA of the matrix P) by a binary counter Cou (serial adder). If the input operands are shown in the additional code (with the sign bits), then all the matrices are formed on the basis of a given bit width of the grid. As an example, Fig. 5 shows the shape matrix and distribution of matrix bits for 24-bit operands A, B, C, D, for diapason (2 - 0).

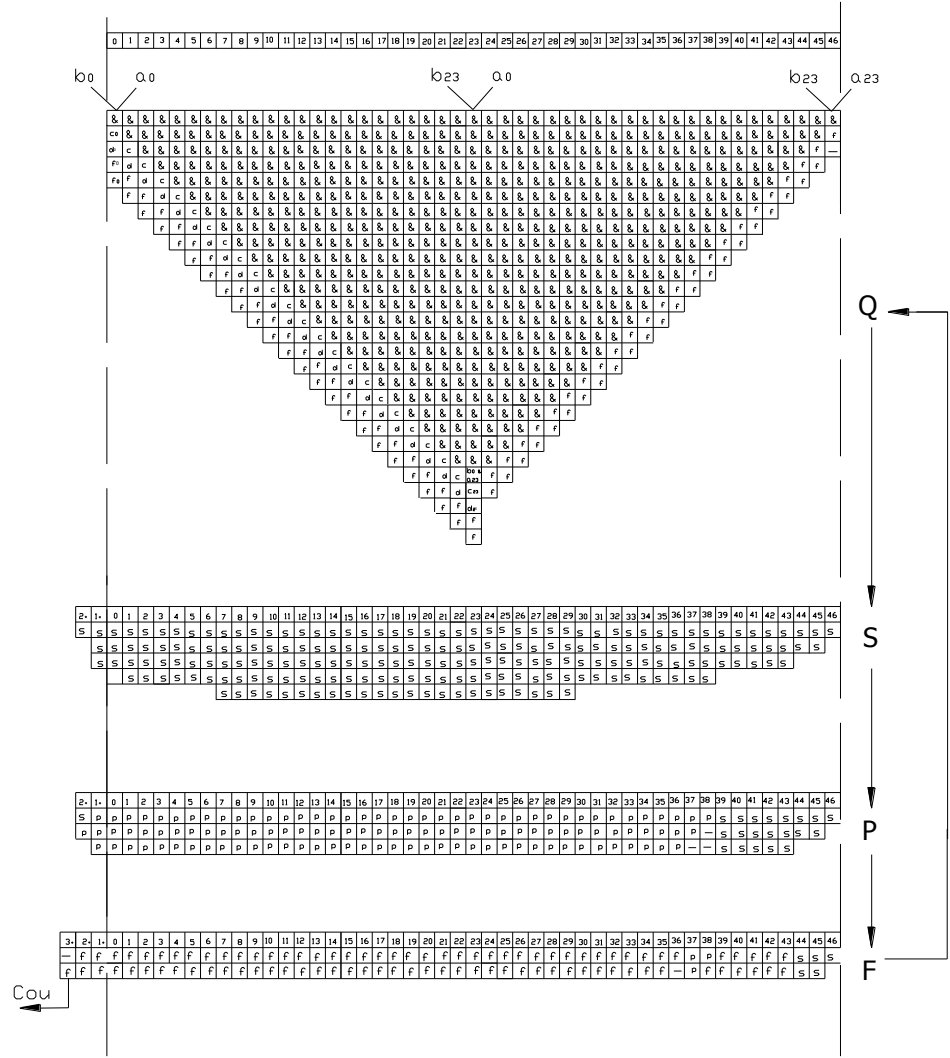

Fig. 5

Computation time operand F (6.1.) in this case is:

$$t_f = t_\& + t_q + t_s + t_p = 14 t_\&,$$

where $t_\&$ - time the formation of the matrix of partial products **&** equal to the delay in one element AND ( two-input);

$t_q = 6 t_\&$ - time convolution matrix **Q** the into matrix **S** (see Table 3.1..);

$t_s = 4 t_\&$ --time convolution matrix **S** the into matrix **P**;

$t_p = 3 t_\&$ - time convolution matrix **P**.

Thus, the operator (6.1) on one iteration for a 24-bit number is calculated during $14 t_\&$. Hardware costs for this structure is approximately 12,500 elements AND (two-input ).

**8. Conclusion.**

The problems of increasing the speed of arithmetic operations on the basis of a distributed arithmetic in the past decades has be given considerable attention. Especially a lot of publication in this area focus on the sphere for the calculation of the difference equations; in digital filters; in structures to solve the equations used in the spectral analysis; in structures for solving matrix equations, etc. [4,5,6]. Most of these works are based on the use in a matrix structure of standard chips. In this paper considered the structure matrix applied to the chips, operating on the basis of a distributed arithmetic. Implementation matrix the structures in of computing devices it involves great hardware costs that hampered their implementation. However, great strides in IC chips allow realization new ideas reviewed in this article and a new approach to the architectures of supercomputers. It should also take into account that proposed in this paper, matrix structures have a greater degree of regularity, which is very important for the integrated technology. In this paper considered the structure matrix allowing multiple increase the speed of arithmetic operations in comparison with classical structures in which the basic element is a parallel adder.

REFERENCES

1. **Евдокимов В.Ф., Стасюк А.И., Щербаков В.И**. Матричные вычислительные устройства. «Наукова думка», 1993. С. 1-152
2. **Нестеренко С.А., Паулин О.Н.** Построение обобщенной модели операции свертки многорядных кодов при цифровой обработке сигналов.//Технология и конструирование в электронной аппаратуре.-2008.-№1-С.20-26
3. **Храпченко В.М.** Методы ускорения арифметических операций, основанные на преобразовании многорядного кода//Вопросы радиоэлектроники. Сер. Электронная вычислительная техника.-1965.-Вып.8-С.121-144
4. **Mintzer,L** «FIR filters with the Xilinx FPGA» FPGA, 92 ACM/SIGDA Workshop on FPGAs hh.129-134
5. **Santoro M.R.** Design and clocking of VLSI multipliers/Stanford Universitey, Computer Systems Laboratory. Report Number: CSL-TR-89-397, October 1989.
6. The Role of Distributed Arithmetic in FPGA-based Signal Processing, Xilinx, Inc.1996